\documentclass{article}
\usepackage{amssymb}

\input{tcilatex}
\begin{document}

\begin{center}
\bigskip \textbf{Regularity properties of nonlinear Schr\"{o}dinger
equations in vector-valued spaces}

\textbf{Veli Shakhmurov}

Department of Mechanical Engineering, Okan University, Akfirat, Tuzla 34959
Istanbul, Turkey,

E-mail: veli.sahmurov@okan.edu.tr

A\textbf{bstract}
\end{center}

In this paper, regularity properties of Cauchy problem for linear and
nonlinear abstract Schr\"{o}dinger equations in vector-valued function
spaces are obtained.


\begin{center}
\ \ \textbf{AMS 2010: 35Q41, 35K15, 47B25, 47Dxx, 46E40 }

\textbf{1. Introduction, definitions}
\end{center}

\bigskip\ Consider the Cauchy problem for nonlinear abstract Schr\"{o}dinger
(NLAS) equations

\begin{equation}
i\partial _{t}u+\Delta u+Au+F\left( u\right) =0,\text{ }x\in R^{n},\text{ }%
t\in \left[ 0,T\right] ,  \tag{1.1}
\end{equation}%
\[
u\left( 0,x\right) =u_{0}\left( x\right) \text{, for a.e. }x\in R^{n} 
\]%
where $A$ is a linear and $F$ is a nonlinear operators in a Banach space $E$%
, $\Delta $ denotes the Laplace operator in $R^{n}$ and $u=$ $u(t,x)$ is the 
$E$-valued unknown function. If $F\left( u\right) =$ $\lambda \left\vert
u\right\vert ^{p}u$ in $\left( 1.1\right) $ we get the nonlinear problem%
\begin{equation}
i\partial _{t}u+\Delta u+Au+\lambda \left\vert u\right\vert ^{p}u=0,\text{ }%
x\in R^{n},\text{ }t\in \left[ 0,T\right] ,  \tag{1.2}
\end{equation}%
\[
u\left( 0,x\right) =u_{0}\left( x\right) \text{, for a.e. }x\in R^{n}, 
\]%
where $p\in \left( 1,\infty \right) $, $\lambda $ is a real number,

\bigskip By rescaling the values of $u$\ it is possible to restrict
attention to the cases $\lambda =1$ or $\lambda =-1$ these call as the
focusing and defocusing abstract Schr\"{o}dinger equations, respectively.
The problem $\left( 1.1\right) $ also contain two critical case. These are
the mass-critical abstract Schr\"{o}dinger equation,%
\[
i\partial _{t}u+\Delta u+Au+\lambda \left\vert u\right\vert ^{\frac{4}{n}}=0,%
\text{ }x\in R^{n},\text{ }t\in \left[ 0,T\right] , 
\]%
which is associated with the conservation of mass,

\[
M\left( u\left( t\right) \right) :=\dint\limits_{R^{n}}\left\Vert u\left(
t,x\right) \right\Vert _{E}^{2}dx 
\]%
and the energy-critical abstract Schr\"{o}dinger equation (in dimensions $%
n>2 $),

\begin{equation}
i\partial _{t}u+\Delta u+Au+\lambda \left\vert u\right\vert ^{\frac{4}{n-2}%
}=0,\text{ }x\in R^{n},\text{ }t\in \left[ 0,T\right] ,  \tag{1.3}
\end{equation}%
which is associated with the conservation of energy,%
\[
E\left( u\left( t\right) \right) :=\dint\limits_{R^{n}}\left[ \frac{1}{2}%
\left\Vert \nabla u\left( t,x\right) \right\Vert _{E}^{2}+\lambda \left( 
\frac{1}{2}-\frac{1}{n}\right) \left\Vert u\left( t,x\right) \right\Vert
_{E}^{\frac{2n}{n-2}}\right] dx. 
\]

Let $\mathbb{N}$ and $\mathbb{C}$\ denote the sets of all natural and
complex numbers, respectively. For $E=\mathbb{C}$ and $A=0$ the problem $%
\left( 1.2\right) $\ become the classical Cauchy problem for nonlinear Schr%
\"{o}dinger (NLS) equations%
\begin{equation}
i\partial _{t}u+\Delta u+\lambda \left\vert u\right\vert ^{p-1}u=0,\text{ }%
x\in R^{n},\text{ }t\in \left[ 0,T\right] ,  \tag{1.4}
\end{equation}%
\[
u\left( 0,x\right) =u_{0}\left( x\right) \text{, for a.e. }x\in R^{n}. 
\]

The existence of solutions and regularity properties of Cauchy problem for
NLS equations studied e.g in $\left[ 5-7\right] $, $\left[ 9-10\right] $, $%
\left[ 15-16\right] $, $\left[ 26\right] $, $\left[ 28\right] $ and the
referances therein.\ In contrast, to the mentioned above results we will
study the regularity properties of the abstract Cauchy problem $\left(
1.1\right) $. Abstract differential equations studied e.g. in $\left[ \text{1%
}\right] $, $\left[ 8\right] $, $\left[ 11\right] $, $\left[ 17\right] $, $%
\left[ 19\right] $, $\left[ 22-25\right] $ and $\left[ 29\right] .$ Since
the Banach space $E$ is arbitrary and $A$ is a possible linear operator, by
choosing $E$ and $A$ we can obtain numerous classes of Schr\"{o}dinger type
equations and its systems which occur in a wide variety of physical systems.
Our main goal is to obtain the exsistence, uniquness and estimates of
solution to the problem $\left( 1.1\right) .$

If we\ choose $E$ a concrete space, for example $E=L^{2}\left( \Omega
\right) $, $A=L,$ where $\Omega $ is a domin in $R^{m}$ with sufficientli
smooth boundary and $L$ is an elliptic operator in $L^{2}\left( \Omega
\right) $ in $\left( 1.2\right) ,$ then we obtain exsistence, uniquness and
the regularity properties of the mixed problem for linear Scredinger equation

\begin{equation}
i\partial _{t}u+\Delta u+Lu=F\left( t,x\right) ,\text{ }t\in \left[ 0,T%
\right] \text{, }x\in R^{n},\text{ }y\in \Omega ,  \tag{1.5}
\end{equation}%
and the followinng NLS equation%
\[
i\partial _{t}u+\Delta u+Lu+\lambda \left\vert u\right\vert ^{p-1}u=0,\text{ 
}t\in \left[ 0,T\right] \text{, }x\in R^{n},\text{ }y\in \Omega ,
\]%
where $u=u\left( t,x,y\right) .$

\ Moreover, let we choose $E=L^{2}\left( 0,1\right) $ and $A$ to be
differential operator with generalized Wentzell-Robin boundary condition
defined by 
\begin{equation}
D\left( A\right) =\left\{ u\in W^{2,2}\left( 0,1\right) ,\text{ }%
B_{j}u=Au\left( j\right) +\dsum\limits_{i=0}^{1}\alpha _{ij}u^{\left(
i\right) }\left( j\right) ,\text{ }j=0,1\right\} ,\text{ }  \tag{1.6}
\end{equation}%
\[
\text{ }Au=au^{\left( 2\right) }+bu^{\left( 1\right) }+cu,
\]%
where $\alpha _{ij}$ are complex numbers, $a=a\left( y\right) ,$ $b=b\left(
y\right) $, $c=c\left( y\right) $ are complex-valued functions. Then, from
the main our theorem, we get the exsistence, uniquness and regularity
properties of Wentzell-Robin type mixed problem for the for linear
Scredinger equation 
\begin{equation}
i\partial _{t}u+\Delta u+a\frac{\partial ^{2}u}{\partial y^{2}}+b\frac{%
\partial u}{\partial y}+cu=F\left( t,x\right) ,\text{ }  \tag{1.7}
\end{equation}%
\ \ \ 

\begin{equation}
B_{j}u=Au\left( t,x,j\right) +\dsum\limits_{i=0}^{1}\alpha _{ij}u^{\left(
i\right) }\left( t,x,j\right) =0\text{, }j=0,1.  \tag{1.8}
\end{equation}%
\begin{equation}
u\left( 0,x,y\right) =u_{0}\left( x,y\right) \text{, for a.e. }x\in R^{n},%
\text{ }y\in \Omega ,  \tag{1.9}
\end{equation}%
and for the following NLS equation,%
\begin{equation}
i\partial _{t}u+\Delta u+a\frac{\partial ^{2}u}{\partial y^{2}}+b\frac{%
\partial u}{\partial y}+cu+F\left( u\right) =0,\text{ }  \tag{1.10}
\end{equation}%
where 
\[
u=u\left( t,x,y\right) \text{, }t\in \left[ 0,T\right] \text{, }x\in R^{n},%
\text{ }y\in \left( 0,1\right) .
\]%
Note that, the regularity properties of Wentzell-Robin type BVP for elliptic
equations were studied e.g. in $\left[ \text{13, 14 }\right] $ and the
references therein. Moreover, if put $E=l_{2}$ and choose $A$ as a infinite
matrix $\left[ a_{mj}\right] $, $m,j=1,2,...,\infty ,$ then from our results
we obtain the exsistence, uniquness and regularity properties of of Cauchy
problem for the linear system of Scredinger equation Consider at first, the
Cauchy problem for infinity many system of linear Schredinger equations 
\begin{equation}
i\partial _{t}u_{m}+\Delta
u_{m}+\sum\limits_{j=1}^{N}a_{mj}u_{j}=F_{j}\left( t,x\right) ,\text{ }t\in %
\left[ 0,T\right] \text{, }x\in R^{n},  \tag{1.11}
\end{equation}%
\[
u_{m}\left( 0,x\right) =u_{m0}\left( x\right) \text{, for a.e. }x\in R^{n},%
\text{ }
\]%
and infinity many system of NLS equation 
\begin{equation}
i\partial _{t}u_{m}+\Delta u_{m}+\sum\limits_{j=1}^{N}a_{mj}u_{j}+F\left(
u_{1},u_{2},...u_{N}\right) =0,\text{ }t\in \left[ 0,T\right] \text{, }x\in
R^{n},  \tag{1.12}
\end{equation}%
\[
u_{m}\left( 0,x\right) =u_{m0}\left( x\right) \text{, for a.e. }x\in R^{n},%
\text{ }
\]%
where $a_{mj}$ are complex numbers, $u_{j}=u_{j}\left( t,x\right) .$

Let $E$ be a Banach space. $L^{p}\left( \Omega ;E\right) $ denotes the space
of strongly measurable $E$-valued functions that are defined on the
measurable subset $\Omega \subset R^{n}$ with the norm

\[
\left\Vert f\right\Vert _{p}=\left\Vert f\right\Vert _{L^{p}\left( \Omega
;E\right) }=\left( \int\limits_{\Omega }\left\Vert f\left( x\right)
\right\Vert _{E}^{p}dx\right) ^{\frac{1}{p}},\text{ }1\leq p<\infty \ . 
\]

For $p=2$ and $E=H$, where $H$ is a Hilbert space, then $L^{p}\left( \Omega
;E\right) $ become $L^{2}\left( \Omega ;H\right) -$ the Hilbert space of $H$%
-valued functions with inner product:%
\[
\left( f,g\right) _{L^{2}\left( \Omega ;H\right) }=\int\limits_{\Omega
}\left( f\left( x\right) ,g\left( x\right) \right) _{H}dx\text{, for any }%
f,g\in L^{2}\left( \Omega ;H\right) . 
\]

Let $L_{t}^{q}L_{x}^{r}\left( E\right) =L_{t}^{q}L_{x}^{r}\left( \left(
a,b\right) \times \Omega ;E\right) $ denotes the space of strongly
measurable $E$-valued functions that are defined on the measurable set $%
\left( a,b\right) \times \Omega $ with the norm 
\[
\left\Vert f\right\Vert _{L_{t}^{q}L_{x}^{r}\left( \left( a,b\right) \times
\Omega ;E\right) }=\left( \dint\limits_{a}^{b}\left[ \int\limits_{\Omega
}\left\Vert f\left( t,x\right) \right\Vert _{E}^{r}dx\right] ^{\frac{q}{r}%
}dt\right) ^{\frac{1}{q}},\text{ }1\leq q,r<\infty \ . 
\]

Let $C\left( \Omega ;E\right) $ denote the space of $E-$valued, bounded
strongly continious functions on $\Omega $ with norm 
\[
\left\Vert u\right\Vert _{C\left( \Omega ;E\right) }=\sup\limits_{x\in
\Omega }\left\Vert u\left( x\right) \right\Vert _{E}. 
\]

$C^{m}\left( \Omega ;E\right) $\ will denote the spaces of $E$-valued
bounded strongly continuous and $m$-times continuously differentiable
functions on $\Omega $ with norm 
\[
\left\Vert u\right\Vert _{C^{m}\left( \Omega ;E\right) }=\max\limits_{0\leq
\left\vert \alpha \right\vert \leq m}\sup\limits_{x\in \Omega }\left\Vert
D^{\alpha }u\left( x\right) \right\Vert _{E}. 
\]

Here, 
\[
O_{R}=\left\{ x\in R^{n},\text{ }\left\vert x\right\vert <R\right\} \text{
for }R>0. 
\]

Let $E_{1}$ and $E_{2}$ be two Banach spaces. $B\left( E_{1},E_{2}\right) $
will denote the space of all bounded linear operators from $E_{1}$ to $%
E_{2}. $ For $E_{1}=E_{2}=E$ it will be denoted by $B\left( E\right) .$

A closed densely defined linear operator\ $A$ is said to be absolute
positive in a Banach\ space $E$ if every real $\lambda ,$ $\left\vert
\lambda \right\vert >\omega $ is in the resolvent set $\rho \left( A\right) $
and for such $\lambda $%
\[
\left\Vert R\left( \lambda ,A\right) \right\Vert _{B\left( E\right) }\leq
M\left( \left\vert \lambda \right\vert -\omega \right) ^{-1}. 
\]

\textbf{Remark 1.1. }It is known\ that if the operator $A$ is absolute
positive in a Banach\ space $E$ and $0\leq \alpha <1$ then it is an
infinitesimal generator of group of bounded linear operator $U_{A}\left(
t\right) $ satisfying 
\[
\left\Vert U_{A}\left( t\right) \right\Vert _{B\left( E\right) }\leq
Me^{\omega \left\vert t\right\vert },\text{ }t\in \left( -\infty ,\infty
\right) , 
\]%
\begin{equation}
\left\Vert A^{\alpha }U_{A}\left( t\right) \right\Vert _{B\left( E\right)
}\leq M\left\vert t\right\vert ^{-\alpha },\text{ }t\in \left( -\infty
,\infty \right)  \tag{1.10}
\end{equation}%
(see e.g. $\left[ \text{19}\right] ,$ $\left[ \text{\S\ 1.6}\right] $,
Theorem 6.3).

Let $E$ be a Banach space. $S=S(R^{n};E)$ denotes $E$-valed Schwartz class,
i.e. the space of all $E$ -valued rapidly decreasing smooth functions on $%
R^{n}$ equipped with its usual topology generated by seminorms. $S(R^{n};%
\mathbb{C})$ denoted by just $S$.

Let $S^{\prime }(R^{n};E)$ denote the space of all continuous linear
operators, $L:S\rightarrow E$, equipped with the bounded convergence
topology. Recall $S(R^{n};E)$ is norm dense in $L^{p}(R^{n};E)$ when $%
1<p<\infty .$

The Banach space\ $E$\ is called a UMD-space and written as $E\in $ UMD if
only if the Hilbert operator 
\[
\left( Hf\right) \left( x\right) =\lim\limits_{\varepsilon \rightarrow
0}\int\limits_{\left\vert x-y\right\vert >\varepsilon }\frac{f\left(
y\right) }{x-y}dy 
\]%
is bounded in the space $L_{p}\left( R,E\right) ,$ $p\in \left( 1,\infty
\right) $ (see e.g. $\left[ 3\right] $, $\left[ 4\right] $). UMD spaces
include $L_{p}$, $l_{p}$ spaces, Lorentz spaces $L_{pq},$ $p,$ $q\in \left(
1,\infty \right) $

Let $F$ denotes the Fourier trasformation, $\hat{u}=Fu$ and 
\[
s\in \mathbb{R},\text{ }\xi =\left( \xi _{1},\xi _{2},...,\xi _{n}\right)
\in R^{n},\text{ }\left\vert \xi \right\vert ^{2}=\dsum\limits_{k=1}^{n}\xi
_{k}^{2}, 
\]%
\[
\langle \xi \rangle :=\left( 1+\left\vert \xi \right\vert ^{2}\right) ^{%
\frac{1}{2}}. 
\]%
Consider $E-$valued Sobolev space $W^{s,p}\left( R^{n};E\right) $ and
homogeneous Sobolev spaces $\mathring{W}^{s,p}\left( R^{n};E\right) $\
defined by respectively, 
\[
W^{s,p}\left( R^{n};E\right) =\left\{ u:u\in S^{\prime }(R^{n};E),\right. 
\text{ } 
\]%
\[
\left\Vert u\right\Vert _{W^{s,p}\left( R^{n};E\right) }=\left. \left\Vert
F^{-1}\left( 1+\left\vert \xi \right\vert ^{2}\right) ^{\frac{s}{2}}\hat{u}%
\right\Vert _{L^{p}\left( R^{n};H\right) }<\infty \right\} , 
\]

\[
\mathring{W}^{s,p}\left( R^{n};E\right) =\left\{ u:u\in S^{\prime
}(R^{n};E),\left\Vert u\right\Vert _{\mathring{W}^{s,p}\left( R^{n};E\right)
}=\left\Vert F^{-1}\left\vert \xi \right\vert ^{s}\right\Vert _{L^{p}\left(
R^{n};H\right) }<\infty \right\} .\text{ } 
\]%
For $\Omega =R^{n}\times G,$ $\mathbf{p=}\left( p_{1},\text{ }p_{2}\right) ,$
$s\in \mathbb{R}$ and $l\in \mathbb{N}$ we define the $E$-valud anisotropic
Sobolev space $W^{s,l,p}\left( \Omega ;E\right) $ by 
\[
W^{s,l,\mathbf{p}}\left( \Omega ;E\right) :=\left\{ u\in S^{\prime }(\Omega
;E),\text{ }\left\Vert u\right\Vert _{W^{s,l,\mathbf{p}}\left( \Omega
\right) }=\left\Vert u\right\Vert _{W^{s,\mathbf{p}}\left( \Omega \right)
}+\left\Vert u\right\Vert _{W^{l,\mathbf{p}}\left( \Omega \right) }\right\}
, 
\]%
where 
\[
\left\Vert u\right\Vert _{W^{s,\mathbf{p}}\left( \Omega ;E\right)
}=\left\Vert F^{-1}\left( 1+\left\vert \xi \right\vert ^{2}\right) ^{\frac{s%
}{2}}\hat{u}\right\Vert _{L^{\mathbf{p}}\left( \Omega ;E\right) }<\infty , 
\]%
\[
\left\Vert u\right\Vert _{W^{l,\mathbf{p}}\left( \Omega ;E\right)
}=\left\Vert u\right\Vert _{L^{\mathbf{p}}\left( \Omega ;E\right)
}+\dsum\limits_{\left\vert \beta \right\vert =l}\left\Vert D_{y}^{\beta
}u\right\Vert _{L^{\mathbf{p}}\left( \Omega ;E\right) }. 
\]%
The similar way, we define homogeneous anisotropic Sobolev spaces $\mathring{%
W}^{s,l,\mathbf{p}}\left( \Omega ;E\right) $ as:%
\[
\mathring{W}^{s,l,\mathbf{p}}\left( \Omega ;E\right) :=\left\{ u\in
S^{\prime }(\Omega ;E),\text{ }\left\Vert u\right\Vert _{W^{s,l,\mathbf{p}%
}\left( \Omega ;E\right) }=\left\Vert u\right\Vert _{W^{s,\mathbf{p}}\left(
\Omega ;E\right) }+\left\Vert u\right\Vert _{\mathring{W}^{l,\mathbf{p}%
}\left( \Omega ;E\right) }\right\} , 
\]%
where 
\[
\left\Vert u\right\Vert _{\mathring{W}^{s,\mathbf{p}}\left( \Omega ;E\right)
}=\left\Vert F^{-1}\left\vert \xi \right\vert ^{s}\hat{u}\right\Vert _{L^{%
\mathbf{p}}\left( \Omega ;E\right) }<\infty . 
\]

Let $A$ be a linear operator in a Banach space $E.$ Consider Sobolev-Lions
type homogeneous and in homogeneous abstract spaces, respectively 
\[
\mathring{W}^{s,p}\left( R^{n};E\left( A\right) ,E\right) =\mathring{W}%
^{s,p}\left( R^{n};E\right) \cap L^{p}\left( R^{n};E\left( A\right) \right) ,%
\text{ } 
\]%
\[
\left\Vert u\right\Vert _{\mathring{W}^{s,p}\left( R^{n};E\left( A\right)
,E\right) }=\left\Vert u\right\Vert _{\mathring{W}^{s,p}\left(
R^{n};E\right) }+\left\Vert u\right\Vert _{L^{p}\left( R^{n};E\left(
A\right) \right) }<\infty , 
\]

\[
W^{s,p}\left( R^{n};E\left( A\right) ,E\right) =W^{s,p}\left( R^{n};E\right)
\cap L^{p}\left( R^{n};E\left( A\right) \right) ,\text{ } 
\]%
\[
\left\Vert u\right\Vert _{W^{s,p}\left( R^{n};E\left( A\right) ,E\right)
}=\left\Vert u\right\Vert _{W^{s,p}\left( R^{n};E\right) }+\left\Vert
u\right\Vert _{L^{p}\left( R^{n};E\left( A\right) \right) }<\infty . 
\]

Sometimes we use one and the same symbol $C$ without distinction in order to
denote positive constants which may differ from each other even in a single
context. When we want to specify the dependence of such a constant on a
parameter, say $\alpha $, we write $C_{\alpha }$.

\textbf{Definition 1.1.} Consider the initial value problem $(1.1)$ for $%
u_{0}\in \mathring{W}^{2,s}\left( R^{n};E\right) $. This problem is critical
when $s=s_{c}:=\frac{n}{2}-\frac{2}{p},$ subcritical when $s>s_{c}$, and
supercritical when $s<s_{c}$.

\textbf{Definition 1.2.} (Solution). A function $\left[ 0,T\right] \times
R^{n}\rightarrow E$ is a (strong) solution to $(1.1)$ if it lies in the
class $C_{t}^{0}W_{x}^{2,s}=C_{t}^{0}\left( I;W_{x}^{2,s}\left(
R^{n};E\left( A\right) ,E\right) \right) $ and obeys the Duhamel formula 
\[
u\left( t\right) =U_{\Delta +A}\left( t\right)
u_{0}+\dint\limits_{0}^{t}U_{\Delta +A}\left( t-s\right) F\left( u\left(
s\right) \right) ds\text{ for all }t\in \left( 0,T\right) , 
\]%
where $U_{\Delta +A}\left( t\right) $ is a bounded group in $E$ generated by
operator $i\left( \Delta +A\right) .$

We write $a\lesssim b$ to indicate that $a\leq Cb$ for some constant $C$,
which is permitted to depend on some parameters.\bigskip 

\begin{center}
\textbf{3. Dispersive and Strichartz inequalities} \textbf{for linear
Schrodinger equation}
\end{center}

\ It can be shown that fundamental solution of the free abstract Schrodinger
equation:%
\begin{equation}
i\partial _{t}u+\Delta u+Au=0,\text{ }t\in \left[ 0,T\right] ,\text{ }x\in
R^{n}  \tag{3.1}
\end{equation}%
can be exspressed as 
\begin{equation}
U_{A+\Delta }\left( t\right) \left( x,y\right) =U_{A}\left( t\right)
U_{\Delta }\left( t\right) \left( x,y\right) \text{, }  \tag{3.2}
\end{equation}%
$U_{A}\left( t\right) $ is a group generated by $iA$ and $U_{\Delta }\left(
t\right) \left( x,y\right) =e^{i\Delta t}\left( x,y\right) $ is a
fundamental solution of the free Schrodinger equation:%
\[
i\partial _{t}u+\Delta u=0,\text{ }x\in R^{n},\text{ }t\in \left[ 0,T\right]
, 
\]%
i.e. 
\begin{equation}
U_{\Delta }\left( t\right) \left( x,y\right) =\left( 4\pi it\right) ^{-\frac{%
n}{2}}e^{i\left\vert x-y\right\vert ^{2}\mid 4t}\text{, }t\neq 0,  \tag{3.3}
\end{equation}

\[
U_{\Delta }\left( t\right) f\left( x\right) =\left( 2\pi it\right) ^{-\frac{n%
}{2}}\dint\limits_{R^{n}}e^{\frac{i\left\vert x-y\right\vert ^{2}}{2t}%
}f\left( y\right) dy. 
\]

\textbf{Lemma 3.1.\ }Let $A$ be absolute positive in a Banach\ space $E$ and 
$0\leq \alpha <1$ . Then the following dispersive inequalites hold 
\begin{equation}
\left\Vert A^{\alpha }U_{\Delta +A}\left( t\right) f\right\Vert
_{L_{x}^{p}\left( R^{n}:E\right) }\lesssim t^{-\left[ n\left( \frac{1}{2}-%
\frac{1}{p}\right) +\alpha \right] }\left\Vert f\right\Vert
_{L_{x}^{p^{\prime }}\left( R^{n}:E\right) },  \tag{3.4}
\end{equation}%
\begin{equation}
\left\Vert A^{\alpha }U_{\Delta +A}\left( t-s\right) f\right\Vert
_{L^{\infty }\left( R^{n};E\right) }\lesssim \left\vert t-s\right\vert
^{-\left( \frac{n}{2}+\alpha \right) }\left\Vert f\right\Vert _{L^{1}\left(
R^{n};E\right) }  \tag{3.5}
\end{equation}%
for $t\neq 0,$ $2\leq p\leq \infty ,$ $\frac{1}{p}+\frac{1}{p^{\prime }}=1.$

\textbf{Proof. }By using $\left( 3.3\right) $ and Young's integral
inequality we have 
\begin{equation}
\left\Vert U_{\Delta }\left( t\right) f\right\Vert _{L_{x}^{p}\left(
R^{n}:E\right) }\lesssim \left\vert t\right\vert ^{-n\left( \frac{1}{2}-%
\frac{1}{p}\right) }\left\Vert f\right\Vert _{L_{x}^{p^{\prime }}\left(
R^{n}:E\right) },  \tag{3.6}
\end{equation}%
\[
\left\Vert U_{\Delta }\left( t\right) f\right\Vert _{L_{x}^{\infty }\left(
R^{n}:E\right) }\lesssim \left\vert t\right\vert ^{-\frac{n}{2}}\left\Vert
f\right\Vert _{L_{x}^{1}\left( R^{n}:E\right) }. 
\]%
By $\left( 1.10\right) $ we get 
\[
\left\Vert A^{\alpha }U_{A}\left( t\right) \right\Vert _{B\left( E\right)
}\lesssim \left\vert t\right\vert ^{-\alpha }\text{, }t\neq 0. 
\]%
By using then the properties of $U_{\Delta +A}\left( t\right) =U_{\Delta
}\left( t\right) $ $U_{A}\left( t\right) $, the estimates $\left( 3.7\right) 
$ and $\left( 3.6\right) $ we obtain $\left( 3.4\right) $ and $\left(
3.5\right) .$

\textbf{Condition 3.1. }Assume $n\geq 1,$%
\[
\frac{2}{q}+\frac{n}{r}\leq \frac{n}{2},\text{ }2\leq q,r\leq \infty \text{
\ and }\left( n,\text{ }q,\text{ }r\right) \neq \left( 2,\text{ }2,\text{ }%
\infty \right) . 
\]

\textbf{Remark 3.1. }If $\frac{2}{q}+\frac{n}{r}=\frac{n}{2},$ then $(q,$ $%
r) $ is called sharp admissible, otherwise $(q,$ $r)$ is called nonsharp
admissible. Note in particular that when $n>2$ the endpoint $\left( 2\text{, 
}\frac{2n}{n--2}\right) $ is called sharp admissible.

\bigskip For a space-time slab $\left[ 0,T\right] \times R^{n}$, we define
the $E-$valued Strichartz norm%
\[
\left\Vert u\right\Vert _{S^{0}\left( I;E\right) }=\sup\limits_{\left(
q,r\right) \text{ admissible}}\left\Vert u\right\Vert
_{L_{t}^{q}L_{x}^{r}\left( I\times R^{n};E\right) }, 
\]%
where $S^{0}\left( \left[ 0,T\right] ;E\right) $ is the closure of all $E-$%
valued test functions under this norm and $N^{0}\left( \left[ 0,T\right]
;E\right) $ denotes the dual of $S^{0}\left( \left[ 0,T\right] ;E\right) .$

Assume $H$ is an abstract Hilbert space and $Q$ is a Hilbert space of
function. Suppose for each $t\in \mathbb{R}$ an operator $U\left( t\right) $%
: $Q\rightarrow L^{2}\left( \Omega ;E\right) $ obeys the following estimates:

\begin{equation}
\left\Vert U\left( t\right) f\right\Vert _{L_{x}^{2}\left( \Omega ;H\right)
}\lesssim \left\Vert f\right\Vert _{Q}  \tag{3.7}
\end{equation}%
for all $t,$ $\Omega \subset R^{n}$ and all $f\in Q;$%
\begin{equation}
\left\Vert U\left( s\right) U^{\ast }\left( t\right) g\right\Vert
_{L_{x}^{\infty }\left( \Omega ;H\right) }\lesssim \left\vert t-s\right\vert
^{-\frac{n}{2}}\left\Vert g\right\Vert _{L_{x}^{1}\left( \Omega ;H\right) } 
\tag{3.8}
\end{equation}%
\begin{equation}
\left\Vert U\left( s\right) U^{\ast }\left( t\right) g\right\Vert
_{L_{x}^{\infty }\left( \Omega ;H\right) }\lesssim \left( 1+\left\vert
t-s\right\vert ^{-\frac{n}{2}}\right) \left\Vert g\right\Vert
_{L_{x}^{1}\left( \Omega ;H\right) }  \tag{3.9}
\end{equation}%
for all $t\neq s$ and all $g\in L_{x}^{1}\left( \Omega ;H\right) .$

For proving the main theorem of this ection, we will use the following
bilinear interpolation result (see $\left[ 2\right] $, Section 3.13.5(b)).

\bigskip \textbf{Lemma 3.2. }Assume $A_{0}$, $A_{1},$ $B_{0}$, $B_{1},$ $%
C_{0}$, $C_{1}$ are Banach spaces and $T$ is a bilinear operator bounded
from ($A_{0}\times B_{0}$, $A_{0}\times B_{1},$ $A_{1}\times B_{0}$ ) into ($%
C_{0}$, $C_{1}$, $C_{1}$), respectively. Then whenever $0<\theta _{0},$ $%
\theta _{1}<\theta <1$ are such that $1\leq \frac{1}{p}+\frac{1}{q}$ and $%
\theta =\theta _{0}+$ $\theta _{1}$, the operator is bounded from 
\[
\left( A_{0}\text{, }A_{1}\right) _{\theta _{0}pr}\times \left( B_{0}\text{, 
}B_{1}\right) _{\theta _{1}qr} 
\]%
to $\left( C_{0}\text{, }C_{1}\right) _{\theta r}.$

By following $\left[ \text{15, Theorem 1.2}\right] $ we have:

\textbf{Theorem 3.1. }Assume $U(t)$ obeys $\left( 3.8\right) $ and $\left(
3.9\right) $. Let $U\left( t\right) $ generates absolute positive
infinitesimal generator operator $A$ and $0\leq \alpha <1.$ Then the
following estimates are hold%
\begin{equation}
\left\Vert U\left( t\right) f\right\Vert _{L_{t}^{q}L_{x}^{r}\left( H\right)
}\lesssim \left\Vert f\right\Vert _{Q},  \tag{3.10}
\end{equation}

\begin{equation}
\left\Vert \dint U^{\ast }\left( s\right) F\left( s\right) ds\right\Vert
_{Q}\lesssim \left\Vert F\right\Vert _{L_{t}^{q^{\prime }}L_{x}^{r^{\prime
}}\left( E^{\ast }\right) },  \tag{3.11}
\end{equation}%
\begin{equation}
\dint\limits_{s<t}\left\Vert A^{\alpha }U\left( t\right) U^{\ast }\left(
s\right) F\left( s\right) ds\right\Vert _{L_{t}^{q}L_{x}^{r}\left( H\right)
}\lesssim \left\Vert F\right\Vert _{L_{t}^{\tilde{q}^{\prime }}L_{x}^{\tilde{%
r}^{\prime }}\left( H\right) },  \tag{3.12}
\end{equation}
for all sharp admissible exponent pairs $\left( q,r\right) $, $\left( \tilde{%
q},\tilde{r}\right) .$ Furthermore, if the decay hypothesis is strengthened
to $(3.9)$, then $(3.10)$, $(3.11)$ and $(3.12)$ hold for all admissible $%
\left( q,\text{ }r\right) $, $\left( \tilde{q},\text{ }\tilde{r}\right) .$

\textbf{Proof. The first step: }Consider the nonendpoint case, i.e. $\left(
q,\text{ }r\right) \neq $ $\left( 2,\text{ }\frac{2n}{n-2}\right) $ and will
show firstly, the estimates $\left( 3.10\right) $, $\left( 3.11\right) .$ By
duality, $(3.10)$ is equivalent to $(3.11)$. By the $TT^{\ast }$ method, $%
(3.11)$ is in turn equivalent to the bilinear form estimate 
\begin{equation}
\left\vert \dint \dint \langle \left( A^{\frac{\alpha }{2}}U\left( s\right)
\right) ^{\ast }F\left( s\right) ,\left( A^{\frac{\alpha }{2}}U\left(
t\right) \right) ^{\ast }G\left( t\right) \rangle dsdt\right\vert \lesssim
\left\Vert F\right\Vert _{L_{t}^{q^{\prime }}L_{x}^{r^{\prime }}\left(
H\right) }\left\Vert G\right\Vert _{L_{t}^{q^{\prime }}L_{x}^{r^{\prime
}}\left( H\right) }.  \tag{3.13}
\end{equation}

By symmetry it suffices to show the to the retarded version of $\left(
3.13\right) $%
\begin{equation}
\left\vert T\left( F,G\right) \right\vert \lesssim \left\Vert F\right\Vert
_{L_{t}^{q^{\prime }}L_{x}^{r^{\prime }}\left( E^{\ast }\right) }\left\Vert
G\right\Vert _{L_{t}^{q^{\prime }}L_{x}^{r^{\prime }}\left( E^{\ast }\right)
},  \tag{3.14}
\end{equation}%
where $T\left( F,G\right) $ is the bilinear form defined by 
\[
T\left( F,G\right) =\dint \dint\limits_{s<t}\langle \left( A^{\frac{\alpha }{%
2}}U\left( s\right) \right) ^{\ast }F\left( s\right) ,\left( A^{\frac{\alpha 
}{2}}U\left( t\right) \right) ^{\ast }G\left( t\right) \rangle dsdt 
\]

By real interpolation between the bilinear form of $\left( 3.7\right) $ and
due to estimate $\left( 1.10\right) $ we get 
\[
\left\vert \langle \left( A^{\frac{\alpha }{2}}U\left( s\right) \right)
^{\ast }F\left( s\right) ,\left( A^{\frac{\alpha }{2}}U\left( t\right)
\right) ^{\ast }G\left( t\right) \rangle \right\vert \lesssim \left\Vert
F\left( s\right) \right\Vert _{L_{x}^{2}}\left\Vert G\left( t\right)
\right\Vert _{L_{x}^{2}}. 
\]%
By using the bilinear form of $\left( 3.8\right) $ and $\left( 1.10\right) $
we have 
\begin{equation}
\left\vert \langle \left( A^{\frac{\alpha }{2}}U\left( s\right) \right)
^{\ast }F\left( s\right) ,\left( A^{\frac{\alpha }{2}}U\left( t\right)
\right) ^{\ast }G\left( t\right) \rangle \right\vert \lesssim  \tag{3.15}
\end{equation}%
\[
\left\vert t-s\right\vert ^{-\frac{n}{2}}\left\Vert F\left( s\right)
\right\Vert _{L_{x}^{1}\left( \Omega ;H\right) }\left\Vert G\left( t\right)
\right\Vert _{L_{x}^{1}\left( \Omega ;H\right) }. 
\]%
In a similar way, we obtain 
\begin{equation}
\left\vert \langle \left( A^{\frac{\alpha }{2}}U\left( s\right) \right)
^{\ast }F\left( s\right) ,\left( A^{\frac{\alpha }{2}}U\left( t\right)
\right) ^{\ast }G\left( t\right) \rangle \right\vert \lesssim  \tag{3.16}
\end{equation}%
\[
\left\vert t-s\right\vert ^{--1-\beta \left( r,r\right) }\left\Vert F\left(
s\right) \right\Vert _{L_{x}^{r^{\prime }}\left( \Omega ;H\right)
}\left\Vert G\left( t\right) \right\Vert _{L_{x}^{r^{\prime }}\left( \Omega
;H\right) }, 
\]%
where $\beta (r,\tilde{r})$ is given by 
\begin{equation}
\beta (r,\tilde{r})=\frac{n}{2}-1-\frac{n}{2}\left( \frac{1}{r}-\frac{1}{%
\tilde{r}}\right) .  \tag{3.17}
\end{equation}

It is clear that $\beta (r,r)\leq 0$ when $n>2.$ In the sharp admissible
case we have 
\[
\frac{1}{q}+\frac{1}{q^{\prime }}=-\beta (r,r), 
\]%
and $(3.14)$ follows from $(3.16)$ and the Hardy-Littlewood-Sobolev
inequality ($[20]$) when $q>q^{\prime }.$

If we are assuming the truncated decay $(3.9)$, then $(3.16)$ can be
improved to 
\begin{equation}
\left\vert \langle \left( A^{\frac{\alpha }{2}}U\left( s\right) \right)
^{\ast }F\left( s\right) ,\left( A^{\frac{\alpha }{2}}U\left( t\right)
\right) ^{\ast }G\left( t\right) \rangle \right\vert \lesssim  \tag{3.18}
\end{equation}%
\[
\left( 1+\left\vert t-s\right\vert \right) ^{--1-\beta \left( r,r\right)
}\left\Vert F\left( s\right) \right\Vert _{L_{x}^{r^{\prime }}\left( \Omega
;H\right) }\left\Vert G\left( t\right) \right\Vert _{L_{x}^{r^{\prime
}}\left( \Omega ;H\right) } 
\]%
and now Young's inequality gives $(3.14)$ when 
\[
-\beta (r,r)+\frac{1}{q}>\frac{1}{q^{\prime }}, 
\]%
i.e. $(q,r)$ is nonsharp admissible. This concludes the proof of $\left(
3.10\right) $ and $(3.11)$ for nonendpoint case.

\textbf{The second step; }It remains to prove $\left( 3.10\right) $ and $%
(3.11)$ for the endpoint case, i.e. when%
\[
\left( q,\text{ }r\right) =\left( 2,\text{ }\frac{2n}{n-2}\right) ,n>2. 
\]
It suffices to show $(3.14)$. \ By decomposing $T(F,G)$ dyadically as $%
\dsum\limits_{j}T_{j}(F,G),$ where the summation is over the integers $%
\mathbb{Z}$ and%
\begin{equation}
T_{j}\left( F,G\right) =\dint\limits_{t-2^{j-1}<s\leq t-2^{j}}\langle \left(
A^{\frac{\alpha }{2}}U\left( s\right) \right) ^{\ast }F\left( s\right)
,\left( A^{\frac{\alpha }{2}}U\left( t\right) \right) ^{\ast }G\left(
t\right) \rangle dsdt  \tag{3.19}
\end{equation}%
we see that it suffices to prove the estimate 
\begin{equation}
\dsum\limits_{j}\left\vert T_{j}(F,G)\right\vert \lesssim \left\Vert
F\right\Vert _{L_{t}^{2}L_{x}^{r^{\prime }}\left( H\right) }\left\Vert
G\right\Vert _{L_{t}^{2}L_{x}^{r^{\prime }}\left( H\right) }\text{.} 
\tag{3.20}
\end{equation}%
For this aim, before we will show the following estimate 
\begin{equation}
\left\vert T_{j}(F,G)\right\vert \lesssim 2^{-j\beta \left( a,b\right)
}\left\Vert F\right\Vert _{L_{t}^{2}L_{x}^{a^{\prime }}\left( H\right)
}\left\Vert G\right\Vert _{L_{t}^{2}L_{x}^{b^{\prime }}\left( H\right) } 
\tag{3.21}
\end{equation}%
for all $j\in \mathbb{Z}$ and all $\left( \frac{1}{a},\frac{1}{b}\right) $\
in a neighbourhood of $\left( \frac{1}{r},\frac{1}{r}\right) $. For proving $%
\left( 3.21\right) $ we will use the real interpolation of $H$-valued
Lebesque space and sequence spaces $l_{q}^{s}\left( H\right) $ (see e.g $%
\left[ \text{27}\right] $ \S\ 1.18.2 and 1.18.6). Indeed, by $\left[ \text{%
27, \S\ 1.18.4.}\right] $ we have 
\begin{equation}
\left( L_{t}^{2}L_{x}^{p_{0}}\left( H\right) ,L_{t}^{2}L_{x}^{p_{1}}\left(
H\right) \right) _{\theta ,2}=L_{t}^{2}L_{x}^{p,2}\left( H\right)  \tag{3.22}
\end{equation}%
whenever $p_{0},$ $p_{1}\in \left[ 1,\infty \right] ,$ $p_{0}\neq p_{1}$ and 
$\frac{1}{p}=\frac{1-\theta }{p_{0}}+\frac{\theta }{p_{1}}$ and $\left(
l_{\infty }^{s_{0}}\left( H\right) ,l_{\infty }^{s_{1}}\left( H\right)
\right) _{\theta ,1}=l_{1}^{s}\left( H\right) $ for $s_{0}$, $s_{1}\in 
\mathbb{R}$, $s_{0}\neq s_{1}$ and%
\[
\frac{1}{s}=\frac{1-\theta }{s_{0}}+\frac{\theta }{s_{1}}, 
\]%
where 
\[
l_{q}^{s}\left( H\right) =\left\{ u=\left\{ u_{j}\right\} _{j+1}^{\infty
},u_{j}\in E\text{, }\left\Vert u\right\Vert _{l_{q}^{s}\left( H\right)
}=\left( \dsum\limits_{j=1}^{\infty }2^{jsq}\left\Vert u_{j}\right\Vert
_{H}^{q}\right) ^{\frac{1}{q}}<\infty \right\} . 
\]

By $\left( 3.22\right) $ the estimate $(3.21)$ can be rewritten as 
\begin{equation}
T:L_{t}^{2}L_{x}^{a^{\prime }}\left( H\right) \times
L_{t}^{2}L_{x}^{b^{\prime }}\left( H\right) \rightarrow l_{\infty }^{\beta
\left( a,b\right) },  \tag{3.23}
\end{equation}%
where $T=\left\{ T_{j}\right\} $ is the vector-valued bilinear operator
corresponding to the $T_{j}.$ We apply Lemma 3.2 to $\left( 3.23\right) $
with $r=1$, $p=q=2$ and arbitrary exponents $a_{0},$ $a_{1}$, $b_{0}$, $%
b_{1} $ such that 
\[
\beta \left( a_{0},b_{1}\right) =\beta \left( a_{1},b_{0}\right) \neq \beta
\left( a_{0},b_{0}\right) . 
\]

Using the real interpolation space identities we obtain%
\[
T:L_{t}^{2}L_{x}^{a^{\prime },2}\left( E^{\ast }\right) \times
L_{t}^{2}L_{x}^{b^{\prime },2}\left( E^{\ast }\right) \rightarrow
l_{1}^{\beta \left( a,b\right) } 
\]%
for all $(a,b)$ in a neighbourhood of $(r,r)$. Applying this to $a=b=r$ and
using the fact that $L^{r^{\prime }}\left( H\right) \subset L^{r^{\prime
},2}\left( H\right) $ we obtain%
\[
T:L_{t}^{2}L_{x}^{a^{\prime },2}\left( H\right) \times
L_{t}^{2}L_{x}^{b^{\prime },2}\left( H\right) \rightarrow l_{1}^{0}\left(
H\right) 
\]%
which implies $\left( 3.21\right) .$

Consider the Cauchy problem for forced Schrodinger equation 
\begin{equation}
i\partial _{t}u+\Delta u+Au=F,\text{ }t\in \left[ 0,T\right] ,\text{ }x\in
R^{n},  \tag{3.24}
\end{equation}

\[
u\left( t_{0},x\right) =u_{0}\left( x\right) ,\text{ }t_{0}\in \left[ 0,T%
\right] , 
\]%
where $A$ is a linear operator in a Hilbert space $H.$

We are now ready to state the standard Strichartz estimates:

\textbf{Theorem 3.2}. Assume the Conditions 3.1 is satisfied and suppose $A$
is absolute positive in $H$. Let $\ 0\leq s\leq 1,$ $0\leq \alpha <1,$ $%
u_{0}\in \mathring{W}^{s,2}\left( R^{n};H\left( A^{\alpha }\right) \right) $%
, $F\in N^{0}\left( \left[ 0,T\right] ;\mathring{W}^{s,2}\left(
R^{n};H\right) \right) $ and let $u$ : $\left[ 0,T\right] \times
R^{n}\rightarrow H$ be a solution to $\left( 3.24\right) $. Then%
\begin{equation}
\left\Vert \left\vert \nabla \right\vert ^{s}u\right\Vert _{S^{0}\left( 
\left[ 0,T\right] ;H\right) }+\left\Vert \left\vert \nabla \right\vert
^{s}A^{\alpha }u\right\Vert _{C^{0}\left( \left[ 0,T\right] ;L^{2}\left(
R^{n};H\right) \right) }\lesssim   \tag{3.25}
\end{equation}%
\[
\left\Vert \left\vert \nabla \right\vert ^{s}A^{\alpha }u_{0}\right\Vert
_{L^{2}\left( R^{n}:H\right) }+\left\Vert \left\vert \nabla \right\vert
^{s}F\right\Vert _{N^{0}\left( \left[ 0,T\right] ;H\right) }.
\]

\textbf{Proof. }We will treat by following $\left[ 15\right] $ and $\left[ 16%
\right] .$ By\textbf{\ }the Duhamel formula the solution $\left( 3.24\right) 
$ can be exspressed as 
\[
u\left( t\right) =U_{\Delta +A}\left( t\right)
u_{0}-i\dint\limits_{0}^{t}U_{\Delta +A}\left( t-s\right) F\left( s\right) ds%
\text{ for all }t\text{, }t_{0}\in \left[ 0,T\right] . 
\]

Let $2\leq q,r,\tilde{q},\tilde{r}$ $\leq \infty $ with 
\[
\frac{2}{q}+\frac{n}{r}=\frac{2}{\tilde{q}}+\frac{n}{\tilde{r}}=\frac{n}{2}. 
\]%
If $n=2,$ we also require that $q,$ $\tilde{q}>2.$ Consider first, the
nonendpoint case. The linear operators in $(3.26)$ and $(3.27)$ are adjoints
of one another; thus, by the method of $TT^{\ast }$ both will follow once we
prove 
\begin{equation}
\left\Vert \dint\limits_{s<t}A^{\alpha }U_{\Delta +A}\left( t-s\right)
F\left( s\right) ds\right\Vert _{L_{t}^{q}L_{x}^{r}\left( H\right) }\lesssim
\left\Vert F\right\Vert _{L_{t}^{q^{\prime }}L_{x}^{r^{\prime }}\left(
H\right) }.  \tag{3.26}
\end{equation}

Apply Theorem 3.1 with $Q=L_{x}^{2}\left( R^{n};H\right) =L_{x}^{2}\left(
H\right) .$ The energy estiamate $\left( 3.10\right) $:%
\[
\left\Vert U_{\Delta +A}\left( t\right) f\right\Vert _{L_{x}^{2}\left(
H\right) }\lesssim \left\Vert f\right\Vert _{L_{x}^{2}\left( H\right) } 
\]%
follows from Plancherel's theorem, the untruncated decay estimate%
\[
\left\Vert U_{\Delta }\left( t-s\right) f\right\Vert _{L_{x}^{\infty }\left(
H\right) }\lesssim \left\vert t-s\right\vert ^{-\frac{n}{2}}\left\Vert
f\right\Vert _{L_{x}^{1}\left( H\right) }, 
\]%
from the equality 
\[
U_{\Delta +A}\left( t\right) f=U_{\Delta }\left( t\right) U_{A}\left(
t\right) f 
\]%
and explicit representation of the Schredinger evolution operator%
\[
U_{\Delta }\left( t\right) f\left( x\right) =\left( 2\pi it\right) ^{-\frac{n%
}{2}}\dint\limits_{R^{n}}e^{\frac{i\left\vert x-y\right\vert ^{2}}{2t}%
}f\left( y\right) dy. 
\]

\bigskip Due to properties of the operator $A$, grops $U_{\Delta +A}\left(
t\right) $ and by the dispersive estimate $(3.4)$ we have 
\[
\left\Vert A^{\alpha }\Phi \right\Vert _{E}\lesssim
\dint\limits_{s<t}\left\Vert A^{\alpha }U_{\Delta +A}\left( t-s\right)
ds\right\Vert _{B\left( H\right) }\left\Vert F\left( s\right) \right\Vert
_{H}ds\lesssim 
\]%
\[
\dint\limits_{\mathbb{R}}\left\vert t-s\right\vert ^{-n\left( \frac{1}{2}-%
\frac{1}{p}\right) -\alpha }\left\Vert F\left( s\right) \right\Vert _{H}ds, 
\]%
where 
\[
\Phi =\dint\limits_{s<t}A^{\alpha }U_{\Delta +A}\left( t-s\right) F\left(
s\right) ds. 
\]

Moreover, from above estimate by the Hardy-Littlewood-Sobolev inequality, we
get

\begin{equation}
\left\Vert A^{\alpha }\Phi \right\Vert _{L_{t}^{q}L_{x}^{r}\left(
R^{n+1};H\right) }\lesssim \left\Vert \dint\limits_{\mathbb{R}}\left\vert
t-s\right\vert ^{-n\left( \frac{1}{2}-\frac{1}{p}\right) -\alpha }\left\Vert
F\left( s\right) \right\Vert _{L_{x}^{r^{\prime }}\left( R^{n};H\right)
}ds\right\Vert _{L_{t}^{q}\left( \mathbb{R}\right) }\lesssim  \tag{3.27}
\end{equation}%
\[
\left\Vert F\right\Vert _{L_{t}^{q_{1}}L_{x}^{r^{\prime }}\left( H\right) }, 
\]%
where 
\[
\frac{1}{q_{1}}=\frac{1}{q}+\frac{1}{p}+\frac{1}{2}-\frac{\alpha }{n}. 
\]

The argument just presented also covers $(3.27)$ in the case $q=\tilde{q},r=%
\tilde{r}$. It allows to consider the estimate in dualized form:%
\begin{equation}
\left\vert \dint \dint\limits_{s<t}\langle U_{\Delta +A}\left( t-s\right)
F\left( s\right) ,G\left( t\right) \rangle ds\right\vert \lesssim \left\Vert
F\right\Vert _{L_{t}^{q^{\prime }}L_{x}^{r^{\prime }}\left( H\right)
}\left\Vert G\right\Vert _{L_{t}^{\tilde{q}_{1}}L_{x}^{\tilde{r}^{\prime
}}\left( H\right) }  \tag{3.28}
\end{equation}%
when 
\[
\frac{1}{\tilde{q}_{1}}=\frac{1}{\tilde{q}}+\frac{1}{\tilde{p}}+\frac{1}{2}-%
\frac{\nu }{n}. 
\]%
The case $\tilde{q}=\infty ,$ $\tilde{r}=2$ follows from $(3.27)$, i.e. 
\begin{equation}
K\lesssim \left\Vert \dint\limits_{s<t}U_{\Delta +A}\left( t-s\right)
F\left( s\right) ds\right\Vert _{L_{t}^{\infty }L_{x}^{2}\left( H\right)
}\left\Vert G\right\Vert _{L_{t}^{1}L_{x}^{2}\left( H\right) }\lesssim 
\tag{3.29}
\end{equation}%
\[
\left\Vert F\right\Vert _{L_{t}^{q_{1}}L_{x}^{r^{\prime }}\left( H\right)
}\left\Vert G\right\Vert _{L_{t}^{1}L_{x}^{2}\left( H\right) }, 
\]%
where 
\[
K=\left\vert \dint \dint\limits_{s<t}\langle U_{\Delta +A}\left( t-s\right)
F\left( s\right) ,G\left( t\right) \rangle ds\right\vert . 
\]%
From $\left( 3.29\right) $ we obtain the esimate $\left( 3.28\right) $ when $%
s=0.$ The general case is obtaind by using the same argument.

Now, consider the endpoint case, i.e. $\left( q,r\right) =\left( 2,\frac{2n}{%
n-2}\right) $. It is suffices to show the following estimates%
\begin{equation}
\left\Vert A^{\alpha }U_{\Delta +A}\left( t\right) u_{0}\right\Vert
_{L_{t}^{q}L_{x}^{r}\left( H\right) }\lesssim \left\Vert Au_{0}\right\Vert
_{W^{s,2}\left( R^{n};H\right) },  \tag{3.30}
\end{equation}%
\begin{equation}
\left\Vert A^{\alpha }U_{\Delta +A}\left( t\right) u_{0}\right\Vert
_{C^{0}\left( L_{x}^{2}\left( H\right) \right) }\lesssim \left\Vert
Au_{0}\right\Vert _{W^{s,2}\left( R^{n};H\right) },  \tag{3.31}
\end{equation}

\begin{equation}
\left\Vert \dint\limits_{s<t}A^{\alpha }U_{\Delta +A}\left( t-s\right)
F\left( s\right) ds\right\Vert _{L_{t}^{q}L_{x}^{r}\left( H\right) }\lesssim
\left\Vert F\right\Vert _{L_{t}^{\tilde{q}^{\prime }}L_{x}^{\tilde{r}%
^{\prime }}\left( H\right) },  \tag{3.32}
\end{equation}%
\begin{equation}
\left\Vert \dint\limits_{s<t}A^{\alpha }U_{\Delta +A}\left( t-s\right)
F\left( s\right) ds\right\Vert _{C^{0}L_{x}^{2}\left( H\right) }\lesssim
\left\Vert F\right\Vert _{L_{t}^{q^{\prime }}L_{x}^{\tilde{r}^{\prime
}}\left( H\right) }.  \tag{3.33}
\end{equation}

Indeed, applying Theorem 3.1 for 
\[
Q=L^{2}\left( R^{n};H\right) ,U\left( t\right) =\chi _{\left[ 0,T\right]
}U_{\Delta +A}\left( t\right) 
\]%
with the energy estimate 
\[
\left\Vert U\left( t\right) f\right\Vert _{L^{2}\left( R^{n};H\right)
}\lesssim \left\Vert f\right\Vert _{L^{2}\left( R^{n};H\right) } 
\]%
which follows from Plancherel's theorem, the untruncated decay estimate $%
\left( 3.8\right) $ and by using of Lemma 3.1 we obtain the estimates $%
\left( 3.30\right) $ and $\left( 3.32\right) .$ Let us temporarily replace
the $C_{t}^{0}L_{x}^{2}\left( H\right) $ norm in estimates $\left(
3.30\right) $, $\left( 3.32\right) $ by the $L_{t}^{\infty }L_{x}^{2}\left(
H\right) .$ Then, all of the above the estimates will follow from Theorem
3.1, once we show that $U\left( t\right) $\ obeys the energy estimate $%
\left( 3.7\right) $ and the truncated decay estimate $(3.9)$. The estimate $%
\left( 3.7\right) $ is obtain immediate from Plancherel's theorem, and $%
\left( 3.9\right) $ follows in a similar way as in $\left[ \text{21,
p.223-224}\right] $. To show that the quantity 
\[
GF\left( t\right) =\dint\limits_{s<t}A^{\alpha }U_{\Delta +A}\left(
t-s\right) F\left( s\right) ds 
\]%
is continuous in $L^{2}\left( R^{n};H\right) ,$ we use the the identity 
\[
GF\left( t+\varepsilon \right) =U\left( \varepsilon \right) GF\left(
t\right) +G\left( \chi _{\left[ t,t+\varepsilon \right] }F\right) \left(
t\right) , 
\]%
the continuity of $U\left( \varepsilon \right) $ as an operator on $%
L^{2}\left( R^{n};H\right) $, and the fact that 
\[
\left\Vert \chi _{\left[ t,t+\varepsilon \right] }F\right\Vert _{L_{t}^{%
\tilde{q}^{\prime }}L_{x}^{\tilde{r}^{\prime }}\left( H\right) }\rightarrow 0%
\text{ as }\varepsilon \rightarrow 0. 
\]

From the estimates $\left( 3.30\right) -\left( 3.33\right) $ we obtain $%
\left( 3.25\right) $ for endpoint case.

\begin{center}
\textbf{4. Strichartz type estimates for solution} \textbf{to nonlinear
Schrodinger equation}
\end{center}

\bigskip Consider the initial-value problem 
\begin{equation}
i\partial _{t}u+\Delta u+Au=F\left( u\right) ,\text{ }x\in R^{n},\text{ }%
t\in \left[ 0,T\right] ,  \tag{4.1}
\end{equation}%
\[
u\left( 0,x\right) =u_{0}\left( x\right) \text{, for a.e. }x\in R^{n} 
\]%
for $p\in \left( 1,\infty \right) $, where $A$ is a linear and $F$ is a
nonlinear operator in a Hilbert space $H$, $\lambda $ is a real number, $%
\Delta $ denotes the Laplace operator in $R^{n}$ and $u=$ $u(t,x)$ is the $H$%
-valued unknown function.

\textbf{Condition 4.1.} Assume that the function $F:$ $H\rightarrow H$ is
continuously differentiable and obeys the power type estimates

\begin{equation}
F\left( u\right) =O\left( \left\Vert u\right\Vert ^{1+p}\right) ,\text{ }%
F_{u}\left( u\right) =O\left( \left\Vert u\right\Vert ^{p}\right) ,\text{ } 
\tag{4.2}
\end{equation}

\begin{equation}
F_{u}\left( \upsilon \right) -F_{u}\left( w\right) =O\left( \left\Vert
\upsilon -w\right\Vert ^{\min \left\{ p,1\right\} }+\left\Vert w\right\Vert
^{\max \left\{ 0,p-1\right\} }\right)  \tag{4.3}
\end{equation}%
for some $p>0,$ where $F_{u}\left( u\right) $\ denotes the derivative of
operator function $F$ with respect to $u\in H.$

From $\left( 4.2\right) $ we obtain 
\begin{equation}
\left\Vert F\left( u\right) -F\left( \upsilon \right) \right\Vert \lesssim
\left\Vert u-\upsilon \right\Vert \left( \left\Vert u\right\Vert
^{p}+\left\Vert \upsilon \right\Vert ^{p}\right) .  \tag{4.4}
\end{equation}

\bigskip \textbf{Remark 4.1. }The model example of a nonlinearity obeying
the conditions above is $F(u)=$ $\left\vert u\right\vert ^{p}u$, for which
the critical homogeneous Sobolev space is $\mathring{W}_{x}^{2,s_{c}}\left(
R^{n};H\right) $ with $s_{c}:=\frac{n}{2}-\frac{2}{p}.$

\textbf{Definition 4.1.} A function $F$ : $\left[ 0,T\right] \times
R^{n}\rightarrow H$ is called a (strong) solution to $(4.1)$ if it lies in
the class 
\[
C_{t}^{0}\left( \left[ 0,T\right] ;\mathring{W}_{x}^{2,s}\left(
R^{n};H\left( A\right) \right) \right) \cap L_{t}^{p+2}L_{x}^{\frac{np\left(
p+2\right) }{4}}\left( \left( \left[ 0,T\right] \times R^{n};H\left(
A\right) \right) \right) 
\]%
and obeys the Duhamel formula 
\[
u\left( t\right) =U_{\Delta +A}\left( t\right)
u_{0}+\dint\limits_{0}^{t}U_{\Delta +A}\left( t-s\right) F\left( u\right)
\left( s\right) ds,\text{ for all }t\in \left[ 0,T\right] . 
\]%
We say that $u$ is a global solution if $T=\infty $.

Let $E$ be a Banach space and $B\left( x,\delta \right) $ denotes the boll
in $R^{n}$ centred in $x$ with radius $\delta $ and $M$ denote the $H-$%
valued Hardy-Littlewood type maximal operator that is defined as:%
\[
Mf\left( x\right) =\sup_{\delta >0}\left[ \mu \left( B\left( x,\delta
\right) \right) \right] ^{-1}\dint\limits_{B\left( x,\delta \right)
}\left\Vert f\left( y\right) \right\Vert _{E}dy. 
\]

For proving the main result of this section we need the following:

By following $\left[ \text{20, Ch.1, \S\ 3, Theorem 1}\right] ,$ we obtain
the following result:

\textbf{Proposition 4.1. }Let $f\in L^{p}\left( R^{n};E\right) $ for $%
1<p\leq \infty $. Then $Mf\left( x\right) \in L^{p}\left( R^{n};E\right) $
and 
\[
\left\Vert M\left( f\right) \right\Vert _{L^{p}\left( R^{n};E\right) }\leq
M_{p}\left\Vert f\right\Vert _{L^{p}\left( R^{n};E\right) }. 
\]

\textbf{Proof. }For $E$ $=\mathbb{R}$, the result is obtained from $\left[ 
\text{20, \S 3, Theorem 1}\right] $. The $E-$valued case can be obtained
from the scalar case by applying it to $\tilde{f}\left( x\right) =\left\Vert
f\left( x\right) \right\Vert _{E}.$

A sequence of random variables $\left\{ r_{k}\right\} _{k\geq 0}$ on $\Omega 
$ is called a Rademacher sequence (see e.g.$\left[ 4\right] $) if%
\[
\mathbb{P}(\{r_{k}=1\}=\mathbb{P}(\{r_{k}=-1\}=\frac{1}{2} 
\]%
for $k\geq 0$ and $\left\{ r_{k}\right\} _{k\geq 0}$ are independent. For
instance, one can take $\Omega =(0,1)$ with the Lebesgue measure and 
\[
r_{k}(t)=sign[\sin (2^{k+1}\pi t)]\text{ for }t\in \Omega . 
\]

Let $\eta \in C_{0}^{\infty }\left( \mathbb{R}\right) $\ nonnegative,
supported in $\sigma =\left\{ \frac{1}{2}<\left\vert \xi \right\vert
<2\right\} $ and satisfying 
\[
\dsum\limits_{j=-\infty }^{\infty }\eta \left( 2^{j}\xi \right) \equiv 1. 
\]

Let $l^{p}\left( E\right) $ denotes $E-$valued sequance space (see e.g $%
\left[ \text{27, \S\ 1.18.1.}\right] $). Define Fourier multiplier operators 
\[
Q_{j}f=F^{-1}\eta \left( 2^{-j}\xi \right) \hat{f}. 
\]

From $\left[ \text{18, Proposition 3.2}\right] $ we have the following
Littlewood-Paley type result for $f\in L^{p}\left( R^{n};E\right) :$

\textbf{Proposition 4.2. }Assume $E$ is UMD space, $p\in \left( 1,\infty
\right) $ and $\left\{ r_{j}\right\} _{j\geq 0}$ is a Rademacher sequence.
Then%
\[
\left\Vert f\right\Vert _{L^{p}\left( R^{n};E\right) }\lesssim \left\Vert
\left\{ r_{j}Q_{j}f\right\} _{j\geq 0}\right\Vert _{L^{p}\left(
R^{n};l^{2}\left( E\right) \right) }\lesssim \left\Vert f\right\Vert
_{L^{p}\left( R^{n};E\right) }. 
\]

Consider the vector-valued version of the Fefferman-Stein type maximal
inequality for $E-$valued functions:

\textbf{Proposition 4.3. }Let $E$ be a Banach space, $1<p<\infty $, $1<q\leq
\infty .$ Then there exists a constant $C\left( p,q\right) $ such that for
all $\left\{ f\right\} _{k\geq 0}\in L^{p}\left( R^{n};E\right) $ one has%
\[
\left\Vert \left\{ Mf\right\} _{k\geq 0}\right\Vert _{L^{p}\left(
R^{n};E\right) }\leq C\left( p,q\right) \left\Vert \left\{ f\right\} _{k\geq
0}\right\Vert _{L^{p}\left( R^{n};l_{q}\left( E\right) \right) }. 
\]

\textbf{Proof. }For $q=\infty $ one uses that

\[
\left\Vert Mf_{k}\left( x\right) \right\Vert _{l^{\infty }\left( E\right)
}\leq M\left\Vert f_{k}\left( x\right) \right\Vert _{_{l^{\infty }\left(
E\right) }}\text{, }x\in R^{n}\text{, }k\geq 0 
\]%
and applies the boundedness of $M$ on $L^{p}\left( R^{n}\right) $ to the
function $\tilde{f}\left( x\right) =\left\Vert f_{k}\left( x\right)
\right\Vert _{l^{\infty }\left( E\right) }.$ If $1<q<1$ and $E$ $=\mathbb{R}$%
, the result can be found in $\left[ \text{20}\right] $(Ch.2, \S\ 1, Theorem
1). The $E-$valued case can be obtained from the scalar case by applying it
to%
\[
\left\Vert \left\{ f_{k}\left( x\right) \right\} _{k\geq 0}\right\Vert
_{E}\subset L^{p}\left( R^{n}\right) . 
\]

Then 
\[
\left\Vert D^{\alpha }f\right\Vert _{L^{r}\left( \mathbb{R};E\right)
}\lesssim \left\Vert \dsum\limits_{j=-\infty }^{\infty }2^{j\alpha
}Q_{j}f\right\Vert _{L^{r}\left( \mathbb{R};E\right) }\lesssim \left\Vert
\dsum\limits_{j=-\infty }^{\infty }2^{2j\alpha }\left\Vert Q_{j}f\left(
.\right) \right\Vert _{E}^{2}\right\Vert _{r} 
\]%
for all $f\in W^{\alpha ,r}\left( \mathbb{R};E\right) $ by multiplier
theorem in $L^{r}\left( R^{n};E\right) $ spaces (see e.g.$\left[ 12\right] $%
) and by Proposition 4.2. Moreover, if the right-hand side is finite then $%
D^{\alpha }f\in L^{r}\left( \mathbb{R};E\right) $ in the sense of $E-$valued
distributions. $Q_{j}$ may be realized as a convolution operator $%
Q_{j}f=\psi _{j}\ast f$, where $\psi _{j}\in S\left( \mathbb{R}\right) $ and%
\begin{equation}
\left\vert \psi _{j}\left( x\right) \right\vert +2^{-j}\left\vert \partial
_{x}\psi _{j}\left( x\right) \right\vert \leq C_{N}2^{j}\left(
1+2^{j}\left\vert x\right\vert \right) ^{-N}  \tag{4.5 }
\end{equation}%
for all $N$ uniformly in $j\in \mathbb{Z}$, and%
\begin{equation}
\dint \psi _{j}\left( x\right) dx=0.  \tag{4.6}
\end{equation}

By following $\left[ \text{7, Proposition 3.1 }\right] $ we obtain.

\textbf{Lemma 4.1. }For any\textbf{\ } $g\in W^{\alpha ,r}\left( \mathbb{R}%
;E\right) ,$%
\[
\left\Vert \tilde{Q}_{j}g\left( y\right) -\tilde{Q}_{j}g\left( x\right)
\right\Vert _{E}\leq C\left\{ 
\begin{array}{c}
2^{j}\left\vert x-y\right\vert Mg\left( x\right) \text{ if }\left\vert
x-y\right\vert \leq C2^{-j} \\ 
Mg\left( x\right) +Mg\left( y\right) \text{ for all }x\text{, }y%
\end{array}%
\right. . 
\]

\textbf{Proof. }Construct also $\tilde{\eta}\in C_{0}^{\infty }\left( \sigma
\right) $ satisfying $\tilde{\eta}\eta \equiv \eta .$ Define 
\[
\tilde{Q}_{j}f=F^{-1}\tilde{\eta}\left( 2^{-j}\xi \right) \hat{f} 
\]%
so that the identity operator may be resolved as 
\[
I=\dsum\limits_{j=-\infty }^{\infty }Q_{j}=\dsum\limits_{j=-\infty }^{\infty
}\tilde{Q}_{j}Q_{j}, 
\]%
and $\tilde{Q}_{j}$ is realized by convolution with a Schwartz function $%
\tilde{\psi}_{j}$ satisfying $(4.5)$ and $(4.6)$.

It is clear that 
\[
\left\Vert \tilde{Q}_{j}g\left( y\right) -\tilde{Q}_{j}g\left( x\right)
\right\Vert _{E}\leq \dint \left\vert \psi _{j}\left( y-z\right) -\psi
_{j}\left( x-z\right) \right\vert \left\Vert g\left( z\right) \right\Vert
_{E}dz. 
\]

For any $x,$ we get 
\[
\left\Vert \tilde{Q}_{j}g\left( x\right) \right\Vert _{E}\leq CMg\left(
x\right) 
\]%
because of $\left( 4.5\right) .$ If $\left\vert x-y\right\vert \leq C2^{-j}$
then 
\[
\left\vert \psi _{j}\left( y-z\right) -\psi _{j}\left( x-z\right)
\right\vert \leq C2^{2j}\left\vert x-y\right\vert \left( 1+2^{j}\left\vert
x-z\right\vert \right) ^{-2}, 
\]%
again by $\left( 4.5\right) $. By a standard calculation this implies the
desired estimate (see e.g $\left[ \text{20, p. 62-63}\right] $).

\textbf{Proposition 4.4. }Assume\textbf{\ }$E$ is a UMD space and $F\in
C^{\left( 1\right) }\left( \mathbb{R};E\right) $. Suppose $\alpha \in \left(
0,1\right) ,$ $1<p,$ $q,$ $r<\infty $\ \ and $r^{-1}=p^{-1}+q^{-1}.$ If $%
u\in L^{\infty }\left( \mathbb{R};E\right) ,$ $D^{\alpha }u\in L^{q}\left( 
\mathbb{R};E\right) $ and $F^{\prime }\left( u\right) \in L^{p}\left( 
\mathbb{R};E\right) $, then $D^{\alpha }\left( F\left( u\right) \right) \in
L^{r}\left( \mathbb{R};E\right) $ and 
\[
\left\Vert D^{\alpha }\left( F\left( u\right) \right) \right\Vert
_{L^{r}\left( \mathbb{R};E\right) }\lesssim \left\Vert F^{\prime }\left(
u\right) \right\Vert _{L^{p}\left( \mathbb{R};E\right) }\left\Vert D^{\alpha
}u\right\Vert _{L^{q}\left( \mathbb{R};E\right) }. 
\]

\textbf{Proof. }In view of $\left( 4.6\right) $ we have%
\[
O_{j}F\left( u\right) \left( x\right) =\dint F\left( u\right) \left(
y\right) \psi _{j}\left( x-y\right) dy=\dint \left[ F\left( u\right) \left(
y\right) -F\left( u\right) \left( x\right) \right] \psi _{j}\left(
x-y\right) dy= 
\]%
\begin{equation}
\dint \left[ \dint\limits_{0}^{1}F^{\prime }\left( tu\left( y\right) +\left(
1-t\right) u\left( x\right) \right) dt\right] \left[ u\left( y\right)
-u\left( x\right) \right] \psi _{j}\left( x-y\right) dy.  \tag{4.7}
\end{equation}

By properties of $E-$valued Hardy-Littlewood maximal operator we get 
\[
\left\Vert \dint\limits_{0}^{1}F^{\prime }\left( tu\left( y\right) +\left(
1-t\right) u\left( x\right) \right) dt\right\Vert _{E}\leq 2M\left(
F^{\prime }\left( u\left( x\right) \right) \right) . 
\]

To estimate $(4.7)$ decompose $u=\dsum\limits_{j}Q_{j}u=\dsum\limits_{j}%
\tilde{Q}_{j}Q_{j}u$ to obtain 
\begin{equation}
\left\Vert Q_{j}F\left( u\right) \left( x\right) \right\Vert _{E}\leq 
\tag{4.8}
\end{equation}%
\[
CM\left( F^{\prime }\left( u\left( x\right) \right) \right)
\dsum\limits_{j=-\infty }^{\infty }\dint \left\Vert \tilde{Q}%
_{j}Q_{j}u\left( y\right) -\tilde{Q}_{j}Q_{j}u\left( x\right) \right\Vert
\left\vert \psi _{j}\left( x-y\right) \right\vert dy. 
\]

Break the sum over $j$ into the cases $j<m$ and $j\geq m$. From $\left(
4.8\right) $ we get that 
\[
\dsum\limits_{j<m}\dint \left\Vert \tilde{Q}_{j}Q_{j}u\left( y\right) -%
\tilde{Q}_{j}Q_{j}u\left( x\right) \right\Vert \left\vert \psi _{j}\left(
x-y\right) \right\vert dy\leq 
\]%
\begin{equation}
C\dsum\limits_{j<m}\dint\limits_{\left\vert x-y\right\vert \leq
2^{-j}}2^{j}\left\vert x-y\right\vert M\left( Q_{j}\right) \left( x\right)
2^{m}\left( 1+2^{m}\right) \left\vert x-y\right\vert ^{-3}dy+  \tag{4.9}
\end{equation}%
\[
C\dsum\limits_{j<m}\dint\limits_{\left\vert x-y\right\vert >2^{-j}}\left[
M\left( Q_{j}\right) \left( x\right) +M\left( Q_{j}\right) \left( y\right) %
\right] 2^{m}\left( 1+2^{m}\right) \left\vert x-y\right\vert ^{-3}dy\leq 
\]%
\[
C\dsum\limits_{j<m}2^{j-m}\left[ M\left( Q_{j}\right) \left( x\right)
+M^{2}\left( Q_{j}\right) \left( x\right) \right] \leq
C\dsum\limits_{j<m}2^{j-m}M^{2}\left( Q_{j}\right) \left( x\right) , 
\]%
where $M^{2}=M\circ M.$

Likewise, we get 
\begin{equation}
\dsum\limits_{j\geq m}\dint \left\Vert \tilde{Q}_{j}Q_{j}u\left( y\right) -%
\tilde{Q}_{j}Q_{j}u\left( x\right) \right\Vert \left\vert \psi _{j}\left(
x-y\right) \right\vert dy\leq C\dsum\limits_{j<m}2^{j-m}M^{2}\left(
Q_{j}\right) \left( x\right) .  \tag{4.10}
\end{equation}

Putting $(4.9)$ and $(4.10)$ into $(4.8)$, we have 
\begin{equation}
\left( \dsum\limits_{m=-\infty }^{\infty }2^{m\alpha }\left\Vert
Q_{m}u\left( x\right) \right\Vert _{E}^{2}\right) ^{\frac{1}{2}}\leq
CM\left( F^{\prime }u\left( x\right) \right) \times  \tag{4.11}
\end{equation}%
\[
\left\{ \dsum\limits_{m}2^{m\alpha }\left[ \dsum\limits_{j<m}2^{j-m}M^{2}%
\left( Q_{j}\right) \left( x\right) +\dsum\limits_{j\geq m}M^{2}\left(
Q_{j}\right) \left( x\right) \right] ^{2}\right\} ^{\frac{1}{2}}\leq 
\]%
\[
CM\left( F^{\prime }u\left( x\right) \right) \dsum\limits_{k=-\infty
}^{\infty }2^{-\varepsilon k}\left( \dsum\limits_{j=-\infty }^{\infty
}2^{j\alpha }\left\Vert Q_{j}u\left( x\right) \right\Vert _{E}^{2}\right) ^{%
\frac{1}{2}} 
\]%
by substituting $m=j-k$ after applying Minkowski's inequality, where 
\[
\varepsilon =2\min \left( \alpha \text{, }1-\alpha \right) >0. 
\]

Finally, from $\left( 4.11\right) $ by using Proposition 4.3 we obtain 
\[
\left\Vert D^{\alpha }\left( F\left( u\right) \right) \right\Vert
_{L^{r}\left( \mathbb{R};E\right) }\leq C\left\Vert \left(
\dsum\limits_{m=-\infty }^{\infty }2^{2m\alpha }\left\Vert Q_{m}u\left(
x\right) \right\Vert _{E}^{2}\right) ^{\frac{1}{2}}\right\Vert _{r}\leq 
\]%
\[
C\left\Vert M\left( F^{\prime }u\left( x\right) \right) \left(
\dsum\limits_{j=-\infty }^{\infty }2^{2j\alpha }\left\Vert M^{2}Q_{j}u\left(
x\right) \right\Vert _{E}^{2}\right) ^{\frac{1}{2}}\right\Vert _{r}\leq 
\]%
\[
C\left\Vert M\left( F^{\prime }u\left( x\right) \right) \right\Vert
_{L^{p}\left( \mathbb{R};E\right) }\left\Vert \left( \dsum\limits_{j=-\infty
}^{\infty }2^{j\alpha }\left\Vert Q_{j}u\left( x\right) \right\Vert
_{E}^{2}\right) ^{\frac{1}{2}}\right\Vert _{q}\leq 
\]%
\[
C\left\Vert F^{\prime }\left( u\right) \right\Vert _{L^{p}\left( \mathbb{R}%
;E\right) }\left\Vert D^{\alpha }u\right\Vert _{L^{q}\left( \mathbb{R}%
;E\right) }.
\]

\textbf{Theorem 4.1. }Assume the Cond\i tons 3.1., 4.1 are satisfied and
suppose $A$ is absolute positive in a Banach\ space $E$. Let $0\leq s\leq 1,$
$0\leq \alpha <1,$ $u_{0}\in \mathring{W}^{s,2}\left( R^{n};H\left(
A^{\alpha }\right) \right) $ and $n\geq 1.$ Then there exists $\eta
_{0}=\eta _{0}\left( n\right) >0$ such that if $0<\eta \leq \eta _{0}$ such
that 
\begin{equation}
\left\Vert \left\vert \nabla \right\vert ^{s}U_{\Delta +A}\left( t\right)
A^{\alpha }u_{0}\right\Vert _{L_{t}^{p+2}L_{x}^{\sigma }\left( \left[ 0,T%
\right] \times R^{n};H\right) }\leq \eta ,  \tag{4.5}
\end{equation}%
then here exists a unique solution $u$ to $\left( 4.1\right) $ on $\left[ 0,T%
\right] \times R^{n}.$ Moreover, the following estimates hold

\begin{equation}
\left\Vert \left\vert \nabla \right\vert ^{s}U_{\Delta +A}A^{\alpha
}u\right\Vert _{L_{t}^{p+2}L_{x}^{\sigma }\left( \left[ 0,T\right] \times
R^{n};H\right) }\leq 2\eta ,  \tag{4.6}
\end{equation}%
\begin{equation}
\left\Vert \left\vert \nabla \right\vert ^{s}u\right\Vert _{S^{0}\left( 
\left[ 0,T\right] \times R^{n};H\right) }+\left\Vert A^{\alpha }u\right\Vert
_{C^{0}\left( \left[ 0,T\right] ;\mathring{W}^{s,2}\left( R^{n};H\right)
\right) }\lesssim  \tag{4.7}
\end{equation}%
\[
\left\Vert A^{\alpha }\left\vert \nabla \right\vert ^{s}u_{0}\right\Vert
_{L_{x}^{2}\left( R^{n};H\right) }+\eta ^{1+p}, 
\]%
\begin{equation}
\left\Vert A^{\alpha }u\right\Vert _{S^{0}\left( \left[ 0,T\right] \times
R^{n};H\right) }\lesssim \left\Vert A^{\alpha }u_{0}\right\Vert
_{L_{x}^{2}\left( R^{n};H\right) },  \tag{4.8}
\end{equation}%
where 
\[
r=r\left( p,n\right) =\frac{2n\left( p+2\right) }{2\left( n-2\right) +np}. 
\]

\textbf{Proof. }We apply the standard fixed point argument. More precisely,
using the Strichartz estimates $\left( 3.25\right) $, we will show that the
solution map $u\rightarrow \Phi (u)$ defined by 
\[
\Phi \left( u\right) \left( t\right) :=U_{A}\left( t\right)
u_{0}+\dint\limits_{0}^{t}U_{A}\left( t-s\right) F\left( u\right) \left(
s\right) ds\text{ for all }t\in \left[ 0,T\right] 
\]%
is a contraction on the set $B_{1}\cap B_{2}$ under the metric given by 
\[
d\left( u,\upsilon \right) =\left\Vert u-\upsilon \right\Vert
_{L_{t}^{p+2}L_{x}^{r}\left( \left[ 0,T\right] \times R^{n};H\right) }, 
\]%
where 
\[
B_{1}=\left\{ u\in W^{\infty ,s_{c},2}=L_{t}^{\infty }W_{x}^{s_{c},2}\left( 
\left[ 0,T\right] \times R^{n};H\left( A^{\alpha }\right) \right) :\right. 
\]

\[
\left. \left\Vert u\right\Vert _{W^{\infty ,s_{c},2}}\leq 2\left\Vert
A^{\alpha }u_{0}\right\Vert _{W_{x}^{s_{c},2}\left( R^{n};H\right) }+C\left(
n\right) \left( 2\eta \right) ^{1+p}\right\} , 
\]%
\[
B_{2}=\left\{ u\in W^{p+2,s_{c},r}=L_{t}^{p+2}W_{x}^{s_{c},r}\left( \left[
0,T\right] \times R^{n};H\right) :\right. 
\]%
\[
\left\Vert A^{\alpha }\left\vert \nabla \right\vert ^{s_{c}}u\right\Vert
_{L_{t}^{p+2}L_{x}^{r}\left( \left[ 0,T\right] \times R^{n};H\right) }\leq
2\eta \text{, and } 
\]%
\[
\left. \left\Vert A^{\alpha }u\right\Vert _{L_{t}^{p+2}L_{x}^{r}\left( \left[
0,T\right] \times R^{n};H\right) }\leq 2C\left( n\right) \text{\ }\left\Vert
A^{\alpha }u_{0}\right\Vert _{L_{x}^{2}\left( R^{n};H\right) }\right\} , 
\]%
here $C(n)$ denotes the constant from the Strichartz inequality in $\left(
3.25\right) .$

Note that both $B_{1}$ and $B_{2}$ are closed in this metric. Using the
Strichartz estimate $\left( 3.25\right) $, Proposition 4.4 and Sobolev
embedding in $H-$valued fractional Sobolev spaces $\left[ 23\right] $, we
get that for $u\in B_{1}\cap B_{2}$, 
\[
\left\Vert \Phi \left( u\right) \right\Vert _{L_{t}^{\infty
}W^{s_{c},2}\left( \left[ 0,T\right] \times R^{n};H\right) }\leq \left\Vert
A^{\alpha }u_{0}\right\Vert _{W_{x}^{s_{c},2}\left( R^{n};H\right) }+ 
\]%
\[
C\left( n\right) \left\Vert \langle \nabla \rangle ^{s_{c}}F\left( u\right)
\right\Vert _{_{L_{t}^{\left( p+2\right) /\left( p+1\right)
}L_{x}^{r_{1}}\left( H\right) }}\leq 
\]%
\[
\left\Vert A^{\alpha }u_{0}\right\Vert _{W_{x}^{s_{c},2}\left(
R^{n};H\right) }+C\left( n\right) \left\Vert \langle \nabla \rangle
^{s_{c}}u\right\Vert _{L_{t}^{p+2}L_{x}^{\sigma }\left( H\right) }\left\Vert
u\right\Vert _{L_{t}^{p+2}L_{x}^{np\left( p+2\right) /4}\left( H\right)
}\leq 
\]%
\[
\left\Vert A^{\alpha }u_{0}\right\Vert _{W_{x}^{s_{c},2}\left(
R^{n};H\right) }+C\left( n\right) \left( 2\eta +2C\left( n\right) \left\Vert
A^{\alpha }u_{0}\right\Vert _{L_{x}^{2}\left( R^{n};H\right) }\right)
\left\Vert \left\vert \nabla \right\vert ^{s_{c}}u\right\Vert
_{L_{t}^{p+2}L_{x}^{r}\left( H\right) }\leq 
\]%
\[
\left\Vert A^{\alpha }u_{0}\right\Vert _{W_{x}^{s_{c},2}\left(
R^{n};H\right) }+C\left( n\right) \left( 2\eta +2C\left( n\right) \left\Vert
A^{\alpha }u_{0}\right\Vert _{L_{x}^{2}\left( R^{n};H\right) }\right) \left(
2\eta \right) ^{p}, 
\]%
where 
\[
L_{t}^{q}L_{x}^{r}\left( H\right) =L_{t}^{q}L_{x}^{r}\left( \left[ 0,T\right]
\times R^{n};H\right) \text{, }r_{1}=r_{1}\left( p,n\right) =\frac{2n\left(
p+2\right) }{2\left( n+2\right) +np}. 
\]

Similarly, 
\[
\left\Vert \Phi \left( u\right) \right\Vert _{L_{t}^{p+2}L_{x}^{r}\left(
H\right) }\leq C\left( n\right) \left\Vert A^{\alpha }u_{0}\right\Vert
_{L_{x}^{2}\left( R^{n};H\right) }+C\left( n\right) \left\Vert u\right\Vert
_{L_{t}^{p+2}L_{x}^{r}\left( H\right) }\leq 
\]%
\[
\left\Vert A^{\alpha }u_{0}\right\Vert _{W_{x}^{s_{c},2}\left(
R^{n};H\right) }+2C^{2}\left( n\right) \left\Vert A^{\alpha
}u_{0}\right\Vert _{L_{x}^{2}\left( R^{n};H\right) }\left( 2\eta \right)
^{p}. 
\]

Arguing as above and invoking $\left( 4.5\right) ,$ we obtain 
\[
\left\Vert \left\vert \nabla \right\vert ^{s_{c}}\Phi \left( u\right)
\right\Vert _{L_{t}^{p+2}L_{x}^{r}}\leq \eta +C\left( n\right) \left\Vert
\left\vert \nabla \right\vert ^{s_{c}}\Phi \left( u\right) \right\Vert
_{L_{t}^{\left( p+2\right) /\left( p+1\right) }L_{x}^{r_{1}}}\leq 
\]%
\[
\eta +C\left( n\right) \left( 2\eta \right) ^{1+p}. 
\]

Thus, choosing $\eta _{0}=\eta _{0}\left( n\right) $ sufficiently small, we
see that for $0<$ $\eta \leq \eta _{0}$ the function $\Phi $ maps the set $%
B_{1}\cap B_{2}$ to itself. To see that it is a contraction, we repeat the
computations above and use $(4.4)$ to obtain 
\[
\left\Vert F\left( u\right) -F\left( \upsilon \right) \right\Vert
_{L_{t}^{p+2}L_{x}^{r}}\leq C\left( n\right) \left\Vert F\left( u\right)
-F\left( \upsilon \right) \right\Vert _{L_{t}^{\left( p+2\right) /\left(
p+1\right) }L_{x}^{r_{1}}}\leq 
\]%
\[
C\left( n\right) \left( 2\eta \right) ^{p}\left\Vert u-\left( \upsilon
\right) \right\Vert _{L_{t}^{p+2}L_{x}^{r}}. 
\]

Thus, choosing $\eta _{0}=\eta _{0}\left( n\right) $ small enough, we can
guarantee that is a contraction on the set $B_{1}\cap B_{2}$. By the
contraction mapping theorem, it follows that has a xed point in $B_{1}\cap
B_{2}$. Since $\Phi $ maps into $C_{t}^{0}W_{x}^{s_{c},2}\left( \left[ 0,T%
\right] \times R^{n};H\right) $ we derive that the fixed point of $\Phi $\
is indeed a solution to $(4.1)$.

In view of Definition 4.1, uniqueness follows from uniqueness in the
contraction mapping theorem.

\begin{center}
\textbf{5.The exsistence and uniquness for the system of Schr\"{o}dinger
equation }
\end{center}

Consider at first, the Cauchy problem for linear system of Schredinger
equations 
\begin{equation}
i\partial _{t}u_{m}+\Delta
u_{m}+\sum\limits_{j=1}^{N}a_{mj}u_{j}=F_{j}\left( t,x\right) ,\text{ }t\in %
\left[ 0,T\right] \text{, }x\in R^{n},  \tag{5.1}
\end{equation}%
\[
u_{m}\left( 0,x\right) =u_{m0}\left( x\right) \text{, for a.e. }x\in R^{n},%
\text{ }
\]
where $u=\left( u_{1},u_{2},...,u_{N}\right) ,$ $u_{j}=u_{j}\left(
t,x\right) ,$ $a_{mj}$ are complex numbers. Let $l_{2}=l_{2}\left( N\right) $
and $l_{2}^{s}=l_{2}^{s}\left( N\right) $ (see $\left[ \text{27, \S\ 1.18}%
\right] $). Let $A$ be the operator in $l_{2}\left( N\right) $ defined by%
\[
\text{ }D\left( A\right) =\left\{ u=\left\{ u_{j}\right\} \in l_{2},\text{ }%
\left\Vert Au\right\Vert _{l_{2}\left( N\right) }=\left(
\sum\limits_{m,j=1}^{N}\left\vert a_{mj}u_{j}\right\vert ^{2}\right) ^{\frac{%
1}{2}}<\infty \right\} ,
\]

\[
A=\left[ a_{mj}\right] \text{, }a_{mj}=a_{jm},\text{ }s>0,\text{ }%
m,j=1,2,...,N,\text{ }N\in \mathbb{N}.
\]

\ From Theorem 3.2 we obtain the following result

\bigskip \textbf{Theorem 5.1. }Assume the Conditions 3.1 are hold. Let $%
0\leq s\leq 1,$ $0\leq \alpha <1,$ $u_{0}\in \mathring{W}^{s,2}\left(
R^{n};D\left( A^{\alpha }\right) \right) $, $F\in N^{0}\left( \left[ 0,T%
\right] ;\mathring{W}^{s,2}\left( R^{n};l_{2}\right) \right) $ and $n\geq 1.$
Let $u$ : $\left[ 0,T\right] \times R^{n}\rightarrow l_{2}\left( N\right) $
be a solution to $\left( 5.1\right) $. Then%
\[
\left\Vert \left\vert \nabla \right\vert ^{s}u\right\Vert _{S^{0}\left( 
\left[ 0,T\right] ;l_{2}\right) }+\left\Vert \left\vert \nabla \right\vert
^{s}A^{\alpha }u\right\Vert _{C^{0}\left( \left[ 0,T\right] ;L^{2}\left(
R^{n};l_{2}\right) \right) }\lesssim 
\]%
\[
\left\Vert \left\vert \nabla \right\vert ^{s}A^{\alpha }u_{0}\right\Vert
_{L^{2}\left( R^{n}:l_{2}\right) }+\left\Vert \left\vert \nabla \right\vert
^{s}F\right\Vert _{N^{0}\left( \left[ 0,T\right] ;l_{2}\right) }.
\]

\ \textbf{Proof.} It is easy to see that $A$ is a symmetric operator in $%
l_{2}$ and other conditions of Theorem 3.2 are satisfied. Hence, from Teorem
4.2 we obtain the conculision.

Consider now, the Cauchy problem $\left( 1.10\right) $. We obtain from
Theorem 4.1 the following result

\textbf{Theorem 5.2. }Assume the Conditions 3.1 and 4.1 are hold. Let $0\leq
s\leq 1,$ $0\leq \alpha <1,$ $u_{0}\in \mathring{W}^{s,2}\left(
R^{n};D\left( A^{\alpha }\right) \right) $ and $n\geq 1.$ Then there exists $%
\eta _{0}=\eta _{0}\left( n\right) >0$ such that if $0<\eta \leq \eta _{0}$
such that 
\[
\left\Vert \left\vert \nabla \right\vert ^{s}U_{\Delta +A}\left( t\right)
A^{\alpha }u_{0}\right\Vert _{L_{t}^{p+2}L_{x}^{\sigma }\left( \left[ 0,T%
\right] \times R^{n};l_{2}\right) }\leq \eta ,
\]%
then here exists a unique solution $u$ to $\left( 1.10\right) $ on $\left[
0,T\right] \times R^{n}.$ Moreover, the following estimates hold

\[
\left\Vert \left\vert \nabla \right\vert ^{s}U_{\Delta +A}A^{\alpha
}u\right\Vert _{L_{t}^{p+2}L_{x}^{\sigma }\left( \left[ 0,T\right] \times
R^{n};l_{2}\right) }\leq 2\eta , 
\]%
\[
\left\Vert \left\vert \nabla \right\vert ^{s}u\right\Vert _{S^{0}\left( 
\left[ 0,T\right] \times R^{n};l_{2}\right) }+\left\Vert A^{\alpha
}u\right\Vert _{C^{0}\left( \left[ 0,T\right] ;\mathring{W}^{s,2}\left(
R^{n};l_{2}\right) \right) }\lesssim 
\]%
\[
\left\Vert A^{\alpha }\left\vert \nabla \right\vert ^{s}u_{0}\right\Vert
_{L_{x}^{2}\left( R^{n};l_{2}\right) }+\eta ^{1+p}, 
\]%
\[
\left\Vert A^{\alpha }u\right\Vert _{S^{0}\left( \left[ 0,T\right] \times
R^{n};l_{2}\right) }\lesssim \left\Vert A^{\alpha }u_{0}\right\Vert
_{L_{x}^{2}\left( R^{n};l_{2}\right) }. 
\]%
where 
\[
\sigma =\sigma \left( p,n\right) =\frac{2n\left( p+2\right) }{2\left(
n-2\right) +np}. 
\]

\ \textbf{Proof.} It is easy to see that $A$ is a symmetric operator in $%
l_{2}$ and other conditions of Theorem 4.1 are satisfied. Hence, from Teorem
4.1 we obtain the conculision.

\begin{center}
\textbf{6.The exsistence and uniquness of solution to anisotropic
Schredinger equation}\ \ \ \ \ \ \ \ \ \ \ \ \ \ \ \ \ \ \ \ \ \ \ \ \ \ \ \
\ \ \ \ \ \ \ \ \ \ \ \ \ \ \ \ \ \ \ \ 
\end{center}

Let $\Omega =R^{n}\times G$, $G\subset R^{d},$ $d\geq 2$ is a bounded domain
with $\left( d-1\right) $-dimensional boundary $\partial G$. Consider at
first, the mixed problem for Schredinger equation

\begin{equation}
i\partial _{t}u+\Delta _{x}u+\sum\limits_{\left\vert \alpha \right\vert \leq
2m}a_{\alpha }\left( y\right) D_{y}^{\alpha }u=F\left( t,x\right) ,\text{ } 
\tag{6.1}
\end{equation}%
\[
\text{ }x\in R^{n},\text{ }y\in G,\text{ }t\in \left[ 0,T\right] ,\text{ }%
p\geq 0,
\]

\begin{equation}
B_{j}u=\sum\limits_{\left\vert \beta \right\vert \leq m_{j}}\ b_{j\beta
}\left( y\right) D_{y}^{\beta }u=0\text{, }x\in R^{n},\text{ }y\in \partial
G,\text{ }j=1,2,...,m,  \tag{6.2}
\end{equation}%
\begin{equation}
u\left( 0,x,y\right) =u_{0}\left( x,y\right) \text{ for }x\in R^{n},\text{ }%
y\in G  \tag{6.3}
\end{equation}%
where $u=u\left( t,x,y\right) $ is a solution, $a_{\alpha },$ $b_{j\beta }$
are the complex valued functions, $\lambda =\pm 1,$ $\alpha =\left( \alpha
_{1},\alpha _{2},...,\alpha _{n}\right) $, $\beta =\left( \beta _{1},\beta
_{2},...,\beta _{n}\right) ,$ $\mu _{i}<2m$ and 
\[
D_{x}^{k}=\frac{\partial ^{k}}{\partial x^{k}},\text{ }D_{j}=-i\frac{%
\partial }{\partial y_{j}},\text{ }D_{y}=\left( D_{1,}...,D_{n}\right) ,%
\text{ }y=\left( y_{1},...,y_{n}\right) . 
\]

$\ $

\bigskip Let%
\[
\xi ^{\prime }=\left( \xi _{1},\xi _{2},...,\xi _{n-1}\right) \in R^{n-1},%
\text{ }\alpha ^{\prime }=\left( \alpha _{1},\alpha _{2},...,\alpha
_{n-1}\right) \in Z^{n},\text{ } 
\]%
\[
\text{ }A\left( y_{0},\xi ^{\prime },D_{y}\right) =\sum\limits_{\left\vert
\alpha ^{\prime }\right\vert +j\leq 2m}a_{\alpha ^{\prime }}\left(
y_{0}\right) \xi _{1}^{\alpha _{1}}\xi _{2}^{\alpha _{2}}...\xi
_{n-1}^{\alpha _{n-1}}D_{y}^{j}\text{ for }y_{0}\in \bar{G} 
\]%
\[
B_{j}\left( y_{0},\xi ^{\prime },D_{y}\right) =\sum\limits_{\left\vert \beta
^{\prime }\right\vert +j\leq m_{j}}b_{j\beta ^{\prime }}\left( y_{0}\right)
\xi _{1}^{\beta _{1}}\xi _{2}^{\beta _{2}}...\xi _{n-1}^{\beta
_{n-1}}D_{y}^{j}\text{ for }y_{0}\in \partial G. 
\]

For $\Omega =R^{n}\times G,$ $\mathbf{p=}\left( p_{1},\text{ }p_{2}\right) ,$
$s\in \mathbb{R}$ and $l\in \mathbb{N}$ let $\mathring{W}^{s,l,p}\left(
\Omega \right) =\mathring{W}^{s,l,p}\left( \Omega ;\mathbb{C}\right) .$

From Theorem 3.2 we obtain the following result

\textbf{Theorem 6.1. }Assume the following conditions be satisfied:

\bigskip (1) $G\in C^{2}$, $a_{\alpha }\in C\left( \bar{G}\right) $ for each 
$\left\vert \alpha \right\vert =2m$ and $a_{\alpha }\in L_{\infty }\left(
G\right) $ for each $\left\vert \alpha \right\vert <2m$;

(2) $b_{j\beta }\in C^{2m-m_{j}}\left( \partial G\right) $ for each $j$, $%
\beta $ and $\ m_{j}<2m$, $\sum\limits_{j=1}^{m}b_{j\beta }\left( y^{\prime
}\right) \sigma _{j}\neq 0,$ for $\left\vert \beta \right\vert =m_{j},$ $%
y^{^{\shortmid }}\in \partial G,$ where $\sigma =\left( \sigma _{1},\sigma
_{2},...,\sigma _{n}\right) \in R^{n}$ is a normal to $\partial G$ $;$

(3) for $y\in \bar{G}$, $\xi \in R^{n}$, $\mu \in S\left( \varphi
_{0}\right) $ for $0\leq \varphi _{0}<\pi $, $\left\vert \xi \right\vert
+\left\vert \mu \right\vert \neq 0$ let $\mu +$ $\sum\limits_{\left\vert
\alpha \right\vert =2m}a_{\alpha }\left( y\right) \xi ^{\alpha }\neq 0$;

(4) for each $y_{0}\in \partial G$ local BVP in local coordinates
corresponding to $y_{0}$:%
\[
\mu +A\left( y_{0},\xi ^{\prime },D_{y}\right) \vartheta \left( y\right) =0, 
\]

\[
B_{j}\left( y_{0},\xi ^{\prime },D_{y}\right) \vartheta \left( 0\right)
=h_{j}\text{, }j=1,2,...,m
\]%
has a unique solution $\vartheta \in C_{0}\left( \mathbb{R}_{+}\right) $ for
all $h=\left( h_{1},h_{2},...,h_{n}\right) \in \mathbb{C}^{n}$ and for $\xi
^{\prime }\in R^{n-1};$

(5) Assume the Conditions 3.1 are hold. Let $0\leq s\leq 1,$ $0\leq \alpha
<1,$ $u_{0}\in \mathring{W}^{s,2}\left( R^{n};D\left( A^{\alpha }\right)
\right) $, $F\in N^{0}\left( \left[ 0,T\right] ;\mathring{W}^{s,2}\left(
R^{n};L^{2}\left( G\right) \right) \right) $ and $n\geq 1.$ Let $u$ : $\left[
0,T\right] \times R^{n}\rightarrow L^{2}\left( G\right) $ be a solution to $%
\left( 6.1\right) -\left( 6.3\right) $. Then%
\[
\left\Vert \left\vert \nabla \right\vert ^{s}u\right\Vert _{S^{0}\left( 
\left[ 0,T\right] ;L^{2}\left( G\right) \right) }+\left\Vert \left\vert
\nabla \right\vert ^{s}A^{\alpha }u\right\Vert _{C^{0}\left( \left[ 0,T%
\right] ;L^{2}\left( R^{n};L^{2}\left( G\right) \right) \right) }\lesssim 
\]%
\[
\left\Vert \left\vert \nabla \right\vert ^{s}A^{\alpha }u_{0}\right\Vert
_{L^{2}\left( R^{n}:L^{2}\left( G\right) \right) }+\left\Vert \left\vert
\nabla \right\vert ^{s}F\right\Vert _{N^{0}\left( \left[ 0,T\right]
;L^{2}\left( G\right) \right) }.
\]

\textbf{Proof. }Let us consider the operator $A$ in $H=L^{2}\left( G\right) $
that are defined by 
\[
D\left( A\right) =\left\{ u\in W^{2m,2}\left( G\right) \text{, }B_{j}u=0,%
\text{ }j=1,2,...,m\text{ }\right\} ,\ Au=\sum\limits_{\left\vert \alpha
\right\vert \leq 2m}a_{\alpha }\left( y\right) D_{y}^{\alpha }u\left(
y\right) .
\]

Then the problem $\left( 6.1\right) -\left( 6.3\right) $ can be rewritten as
the problem $\left( 4.1\right) $, where $u\left( x\right) =u\left(
x,.\right) ,$ $f\left( x\right) =f\left( x,.\right) $,\ $x\in R^{n}$ are the
functions with values in\ $H=L^{2}\left( G\right) $. By virtue of $\left[ 
\text{8, Theorem 8.2}\right] ,$ operator $A+\mu $ is absolute positive in $%
L^{2}\left( G\right) $ for sufficiently large $\mu >0$. Moreover, in view of
(1)-(5) all conditons of Theorem 3.2 are hold. Then Theorem 3.2 implies the
assertion.

Consider now, the mixed problem for nonlinear Schrodinger equation

\begin{equation}
i\partial _{t}u+\Delta _{x}u+\sum\limits_{\left\vert \alpha \right\vert \leq
2m}a_{\alpha }\left( y\right) D_{y}^{\alpha }u+\lambda \left\vert
u\right\vert ^{p}u=0,\text{ }  \tag{6.4}
\end{equation}%
\[
\text{ }x\in R^{n},\text{ }y\in G,\text{ }t\in \left[ 0,T\right] ,\text{ }%
p\geq 0,
\]

\begin{equation}
B_{j}u=\sum\limits_{\left\vert \beta \right\vert \leq m_{j}}\ b_{j\beta
}\left( y\right) D_{y}^{\beta }u=0\text{, }x\in R^{n},\text{ }y\in \partial
G,\text{ }j=1,2,...,m,  \tag{6.5}
\end{equation}%
\begin{equation}
u\left( 0,x,y\right) =u_{0}\left( x,y\right) \text{ for }x\in R^{n},\text{ }%
y\in G  \tag{6.6}
\end{equation}

\textbf{Theorem 6.2}. Assume the following conditions be satisfied:

\bigskip (1) $G\in C^{2}$, $a_{\alpha }\in C\left( \bar{G}\right) $ for each 
$\left\vert \alpha \right\vert =2m$ and $a_{\alpha }\in L_{\infty }\left(
G\right) $ for each $\left\vert \alpha \right\vert <2m$;

(2) $b_{j\beta }\in C^{2m-m_{j}}\left( \partial G\right) $ for each $j$, $%
\beta $ and $\ m_{j}<2m$, $\sum\limits_{j=1}^{m}b_{j\beta }\left( y^{\prime
}\right) \sigma _{j}\neq 0,$ for $\left\vert \beta \right\vert =m_{j},$ $%
y^{^{\shortmid }}\in \partial G,$ where $\sigma =\left( \sigma _{1},\sigma
_{2},...,\sigma _{n}\right) \in R^{n}$ is a normal to $\partial G$ $;$

(3) for $y\in \bar{G}$, $\xi \in R^{n}$, $\mu \in S\left( \varphi
_{0}\right) $ for $0\leq \varphi _{0}<\pi $, $\left\vert \xi \right\vert
+\left\vert \mu \right\vert \neq 0$ let $\mu +$ $\sum\limits_{\left\vert
\alpha \right\vert =2m}a_{\alpha }\left( y\right) \xi ^{\alpha }\neq 0$;

(4) for each $y_{0}\in \partial G$ local BVP in local coordinates
corresponding to $y_{0}$:%
\[
\mu +A\left( y_{0},\xi ^{\prime },D_{y}\right) \vartheta \left( y\right) =0, 
\]

\[
B_{j}\left( y_{0},\xi ^{\prime },D_{y}\right) \vartheta \left( 0\right)
=h_{j}\text{, }j=1,2,...,m 
\]%
has a unique solution $\vartheta \in C_{0}\left( \mathbb{R}_{+}\right) $ for
all $h=\left( h_{1},h_{2},...,h_{n}\right) \in \mathbb{C}^{n}$ and for $\xi
^{\prime }\in R^{n-1};$

(5) Assume the Condition 3.1 are hold. Let $0\leq s\leq 1,$ $0\leq \alpha <1,
$ $u_{0}\in \mathring{W}^{s,2,2m}\left( R^{n}\times G\right) $ and $n\geq 1.$

Then there exists $\eta _{0}=\eta _{0}\left( n\right) >0$ such that if $%
0<\eta \leq \eta _{0}$ such that 
\[
\left\Vert \left\vert \nabla \right\vert ^{s}U_{\Delta +A}\left( t\right)
A^{\alpha }u_{0}\right\Vert _{L_{t}^{p+2}L_{x}^{\sigma }L_{y}^{2}\left( 
\left[ 0,T\right] \times R^{n}\times G\right) }\leq \eta ,
\]%
then here exists a unique solution $u$ to $\left( 6.4\right) -\left(
6.6\right) $ on $\left[ 0,T\right] \times R^{n}.$ Moreover, the following
estimates hold

\[
\left\Vert \left\vert \nabla \right\vert ^{s}U_{\Delta +A}A^{\alpha
}u\right\Vert _{L_{t}^{p+2}L_{x}^{\sigma }L_{y}^{2}\left( \left[ 0,T\right]
\times R^{n}\times G\right) }\leq 2\eta , 
\]%
\[
\left\Vert \left\vert \nabla \right\vert ^{s}u\right\Vert
_{S^{0}L_{y}^{2}\left( \left[ 0,T\right] \times R^{n}\times G\right)
}+\left\Vert A^{\alpha }u\right\Vert _{C^{0}\left( \left[ 0,T\right] ;%
\mathring{W}^{s,2}\left( R^{n}\times G\right) \right) }\lesssim 
\]%
\[
\left\Vert A^{\alpha }\left\vert \nabla \right\vert ^{s}u_{0}\right\Vert
_{L_{x,y}^{2}\left( R^{n}\times G\right) }+\eta ^{1+p}, 
\]%
\[
\left\Vert A^{\alpha }u\right\Vert _{S^{0}L_{y}^{2}\left( \left[ 0,T\right]
\times R^{n}\times G\right) }\lesssim \left\Vert A^{\alpha }u_{0}\right\Vert
_{L_{x,y}^{2}\left( R^{n}\times G\right) }. 
\]%
where 
\[
\sigma =\sigma \left( p,n\right) =\frac{2n\left( p+2\right) }{2\left(
n-2\right) +np}. 
\]

\ \textbf{Proof. }The problem $\left( 6.4\right) -\left( 6.6\right) $ can be
rewritten as the problem $\left( 1.1\right) $, where $u\left( x\right)
=u\left( x,.\right) ,$ $f\left( x\right) =f\left( x,.\right) $,\ $x\in R^{n}$
are the functions with values in\ $H=L^{2}\left( G\right) $. By virtue of $%
\left[ \text{8, Theorem 8.2}\right] ,$ operator $A+\mu $ is absolute
positive in $L^{2}\left( G\right) $ for sufficiently large $\mu >0$.
Moreover, in view of (1)-(5) all conditons of Theorem 4.1 are hold. Then
Theorem 4.1 implies the assertion.

\begin{center}
\textbf{7.} \textbf{The Wentzell-Robin type mixed problem for Schr\"{o}%
dinger equations}
\end{center}

Consider at first, the linear problem $\left( 1.7\right) -\left( 1.9\right) $%
. From Theorem 3.2 we obtain the following result

\textbf{Theorem 7.1. } Suppose the the following conditions are satisfied:

(1)\ $a$ is positive, $b$ is a real-valued functions on $\left( 0,1\right) $%
. Moreover$,$ $a\left( .\right) \in C\left( 0,1\right) $ and%
\[
\exp \left( -\dint\limits_{\frac{1}{2}}^{x}b\left( t\right) a^{-1}\left(
t\right) dt\right) \in L_{1}\left( 0,1\right) ; 
\]

(2) Assume the Conditions 3.1 and 4.1 are hold. Let $0\leq s\leq 1,$ $0\leq
\alpha <1,$ $F\in N^{0}\left( \left[ 0,T\right] ;\mathring{W}^{s,2}\left(
R^{n};L^{2}\left( 0,1\right) \right) \right) $ $u_{0}\in \mathring{W}%
^{s,2,2}\left( R^{n}\times \left( 0,1\right) \right) ,$  and $n\geq 1.$

Let $u$ : $\left[ 0,T\right] \times R^{n}\rightarrow L^{2}\left( G\right) $
be a solution to $\left( 1,7\right) -\left( 1.9\right) $. Then%
\[
\left\Vert \left\vert \nabla \right\vert ^{s}u\right\Vert _{S^{0}\left( 
\left[ 0,T\right] ;L^{2}\left( 0,1\right) \right) }+\left\Vert \left\vert
\nabla \right\vert ^{s}A^{\alpha }u\right\Vert _{C^{0}\left( \left[ 0,T%
\right] ;L^{2}\left( R^{n};L^{2}\left( 0,1\right) \right) \right) }\lesssim 
\]%
\[
\left\Vert \left\vert \nabla \right\vert ^{s}A^{\alpha }u_{0}\right\Vert
_{L^{2}\left( R^{n}:L^{2}\left( 0,1\right) \right) }+\left\Vert \left\vert
\nabla \right\vert ^{s}F\right\Vert _{N^{0}\left( \left[ 0,T\right]
;L^{2}\left( 0,1\right) \right) }.
\]

\ \textbf{Proof.} Let $H=L^{2}\left( 0,1\right) $ and $A$ is a operator
defined by $\left( 4.1\right) .$ Then the problem $\left( 1.7\right) -\left(
1.9\right) $ can be rewritten as the problem $\left( 1.2\right) $. By virtue
of $\left[ \text{13, 14}\right] $ the operator $A$ generates analytic
semigroup in $L^{2}\left( 0,1\right) $. Hence, by virtue of (1)-(5) all
conditons of Theorem 3.2 are satisfied. Then Theorem 3.2 implies the
assertion.

Consider now, the problem $\left( 1.7\right) -\left( 1.9\right) $. In this
section, from Theorem 4.1 we obtain the following result:

\bigskip \textbf{Theorem 7.2. } Suppose the the following conditions are
satisfied:

(1)\ $a$ is positive, $b$ is a real-valued functions on $\left( 0,1\right) $%
. Moreover$,$ $a\left( .\right) \in C\left( 0,1\right) $ and%
\[
\exp \left( -\dint\limits_{\frac{1}{2}}^{x}b\left( t\right) a^{-1}\left(
t\right) dt\right) \in L_{1}\left( 0,1\right) ; 
\]

(2) Assume the Conditions 3.1 and 4.1 are hold. Let $0\leq s\leq 1,$ $0\leq
\alpha <1,$ $u_{0}\in \mathring{W}^{s,2,2}\left( R^{n}\times \left(
0,1\right) \right) $ and $n\geq 1.$

Then there exists $\eta _{0}=\eta _{0}\left( n\right) >0$ such that if $%
0<\eta \leq \eta _{0}$ such that 
\[
\left\Vert \left\vert \nabla \right\vert ^{s}U_{\Delta +A}\left( t\right)
A^{\alpha }u_{0}\right\Vert _{L_{t}^{p+2}L_{x}^{\sigma }L_{y}^{2}\left( 
\left[ 0,T\right] \times R^{n}\times \left( 0,1\right) \right) }\leq \eta ,
\]%
then here exists a unique solution $u$ to $\left( 1.8\right) -\left(
1.10\right) $ on $\left[ 0,T\right] \times R^{n}.$ Moreover, the following
estimates hold

\[
\left\Vert \left\vert \nabla \right\vert ^{s}U_{\Delta +A}A^{\alpha
}u\right\Vert _{L_{t}^{p+2}L_{x}^{\sigma }L_{y}^{2}\left( \left[ 0,T\right]
\times R^{n}\times \left( 0,1\right) \right) }\leq 2\eta , 
\]%
\[
\left\Vert \left\vert \nabla \right\vert ^{s}u\right\Vert
_{S^{0}L_{y}^{2}\left( \left[ 0,T\right] \times R^{n}\times \left(
0,1\right) \right) }+\left\Vert A^{\alpha }u\right\Vert _{C^{0}\left( \left[
0,T\right] ;\mathring{W}^{s,2}\left( R^{n}\times 0,1\right) \right)
}\lesssim 
\]%
\[
\left\Vert A^{\alpha }\left\vert \nabla \right\vert ^{s}u_{0}\right\Vert
_{L_{x,y}^{2}\left( R^{n}\times \left( 0,1\right) \right) }+\eta ^{1+p}, 
\]%
\[
\left\Vert A^{\alpha }u\right\Vert _{S^{0}L_{y}^{2}\left( \left[ 0,T\right]
\times R^{n}\times \left( 0,1\right) \right) }\lesssim \left\Vert A^{\alpha
}u_{0}\right\Vert _{L_{x,y}^{2}\left( R^{n}\times \left( 0,1\right) \right)
}. 
\]%
where 
\[
\sigma =\sigma \left( p,n\right) =\frac{2n\left( p+2\right) }{2\left(
n-2\right) +np}. 
\]

\ \textbf{Proof.} Let $H=L^{2}\left( 0,1\right) $ and $A$ is a operator
defined by $\left( 1.4\right) .$ Then the problem $\left( 1.8\right) -\left(
1.10\right) $ can be rewritten as the problem $\left( 4.1\right) $. By
virtue of $\left[ \text{13, 14}\right] $ the operator $A$ generates analytic
semigroup in $L^{2}\left( 0,1\right) $. Hence, by virtue of (1)-(5) all
conditons of Theorem 4.1 are satisfied. Then Theorem 4.1 implies the
assertion.

\textbf{References}\ \ 

\begin{enumerate}
\item H. Amann, Linear and quasi-linear equations,1, Birkhauser, Basel 1995.

\item J. Bergh and J. Lofstrom, Interpolation spaces: An introduction,
Springer-Verlag, New York, 1976.

\item D. L. Burkholder, A geometrical conditions that implies the existence
certain singular integral of Banach space-valued functions, Proc. conf.
Harmonic analysis in honor of Antoni Zygmund, Chicago, 1981,Wads Worth,
Belmont, (1983), 270-286.

\item J. Bourgain, Vector--Valued Singular integrals and the H1--BMO
duality. In: Probability theory and harmonic analysis, pp. 1 -- 19, Pure
Appl. Math. 98, Marcel Dekker, 1986.

\item J. Bourgain, Global solutions of nonlinear Schrodinger equations.
American Mathematical, Society Colloquium Publications, 46. American
Mathematical Society, Providence, RI, 1999, MR1691575.

\item T. Cazenave and F. B. Weissler, The Cauchy problem for the critical
nonlinear Schrodinger equation in $H^{s}$. Nonlinear Anal. 14 (1990),
807-836.

\item M. Christ and M. Weinstein, Dispersion of small amplitude solutions of
the generalized Korteweg-de Vries equation. J. Funct. Anal. 100 (1991),
87-109.

\item R. Denk, M. Hieber, J. Pr\"{u}ss J, $R$-boundedness, Fourier
multipliers and problems of elliptic and parabolic type, Mem. Amer. Math.
Soc. 166 (2003), n.788.

\item J. Ginibre and G. Velo, Smoothing properties and retarded estimates
for some dispersive evolution equations, Comm. Math. Phys. 123 (1989),
535--573.

\item L. Escauriaza, C. E. Kenig, G. Ponce, and L. Vega, Hardy's uncertainty
principle, convexity and Schr\"{o}dinger evolutions, J. European Math. Soc.
10, 4 (2008) 883--907.

\item J. A. Goldstain, Semigroups of linear operators and applications,
Oxford University Press, Oxfard, 1985.

\item M. Girardi and L. Weis. Operator-valued Fourier multiplier theorems on
Lp(X) and geometry of Banach spaces. J. Funct. Anal., 204(2):320--354, 2003.

\item A. Favini, G. R. Goldstein, J. A. Goldstein and S. Romanelli,
Degenerate second order differential operators generating analytic
semigroups in $L_{p}$ and $W^{1,p}$, Math. Nachr. 238 (2002), 78 --102.

\item V. Keyantuo, M. Warma, The wave equation with Wentzell--Robin boundary
conditions on Lp-spaces, J. Differential Equations 229 (2006) 680--697.

\item M. Keel and T. Tao, Endpoint Strichartz estimates. Amer. J. Math. 120
(1998), 955-980.

\item R. Killip and M. Visan, Nonlinear Schrodinger equations at critical
regularity, Clay Mathematics Proceedings, v. 17, 2013.

\item Lunardi A., Analytic Semigroups and optimal regularity in parabolic
problems, Birkhauser, 2003.

\item M. Meyries, M. Veraar, Pointwise multiplication on vector-valued
function spaces with power weights, J. Fourier Anal. Appl. 21 (2015)(1),
95--136.

\item A. Pazy, Semigroups of linear operators and applications to partial
differential equations. Springer, Berlin, 1983.

\item E. M. Stein, Singular Integrals and differentiability properties of
functions, Princeton Univ. Press, Princeton. NJ, 1970.

\item C. D. Sogge, Fourier Integrals in Classical Analysis, Cambridge
University Press, 1993.

\item V. B. Shakhmurov, Nonlinear abstract boundary value problems in
vector-valued function spaces and applications, Nonlinear Anal-Theor., v.
67(3) 2006, 745-762.

\item V. B. Shakhmurov, Embedding and separable differential operators in
Sobolev-Lions type spaces, Mathematical Notes, (84)6 2008, 906-926.

\item R. Shahmurov, On strong solutions of a Robin problem modeling heat
conduction in materials with corroded boundary, Nonlinear Anal., Real World
Appl., v.13, (1), 2011, 441-451.

\item R. Shahmurov, Solution of the Dirichlet and Neumann problems for a
modified Helmholtz equation in Besov spaces on an annuals, J. Differential
equations, v. 249(3), 2010, 526-550.

\item R. S. Strichartz, Restriction of Fourier transform to quadratic
surfaces and deay of solutionsof wave equations. Duke Math. J. 44 (1977),
705\{714. MR0512086A.

\item H. Triebel, Interpolation theory, Function spaces, Differential
operators, North-Holland, Amsterdam, 1978.

\item T. Tao, Nonlinear dispersive equations. Local and global analysis.
CBMS Regional conference series in mathematics, 106. American Mathematical
Society, Providence, RI, 2006., MR2233925

\item S. Yakubov and Ya. Yakubov, Differential-operator Equations. Ordinary
and Partial \ Differential Equations, Chapman and Hall /CRC, Boca Raton,
2000.
\end{enumerate}

\end{document}